\documentclass[letterpaper, 11 pt]{amsart}
\usepackage{Macros,subfig,standalone, amsthm, amsmath, amssymb, amsfonts, comment}

\title{Tableau Stabilization and Lattice paths}
\author{ Connor Ahlbach, Jacob David, Suho Oh, Christopher Wu }
\date{June 2020}
\usepackage[font=footnotesize]{caption}
\usepackage{xcolor}

\newcommand{\newword}[1]{\textbf{\emph{#1}}}
\newcommand{\subseq}{\preceq}
\newtheorem{innercustomthm}{Theorem}
\newenvironment{customthm}[1]
  {\renewcommand\theinnercustomthm{#1}\innercustomthm}
  {\endinnercustomthm}
\usepackage{graphicx}
\begin{document}

\begin{abstract} If one attaches shifted copies of a skew tableau to the right of itself and rectifies, at a certain point the copies no longer experience vertical slides, a phenomenon called tableau stabilization. While tableau stabilization was originally developed to construct the sufficiently large rectangular tableaux fixed by given powers of promotion, the purpose of this paper is to improve the original bound on tableau stabilization to the number of rows of the skew tableau. In order to prove this bound, we encode increasing subsequences as lattice paths and show that various operations on these lattice paths weakly increase the maximum combined length of the increasing subsequences.

\end{abstract}

\maketitle

\section{Introduction}

Tableau stabilization was introduced in \cite{TabStab} in order to construct sufficiently large rectangular tableaux fixed by given powers of promotion. Rhoades first counted the number of rectangular tableaux fixed by the powers of promotion by exhibiting a remarkable cyclic sieving phenomenon \cite{MR2087303} for the action of promotion on rectangular tableaux \cite{MR2557880}. Since Rhoades's work, there has been significant interest in promotion and cyclic sieving for various kinds of tableaux \cite{PromSSYT},\cite{SkewCSP},\cite{BENNETT201462},\cite{Fontaine2012CyclicSR},\cite{Oh2019CrystalsST},\cite{PECHENIK2014357},\cite{10.2307/23513431}.


But while many of these CSPs count fixed points of promotion, they say nothing about what the actual fixed points are. Purbhoo and Rhee first found all rectangular tableaux of shape $ (a^b) $ fixed by $ a $ promotions for $ b \ge a $ \cite{MR3625918}. More recently, Ahlbach exhibited all sufficiently large rectangular tableaux fixed by a given power of promotion by applying the rectification operator to skew tableaux formed by attaching shifted copies of a skew tableaux to itself \cite{TabStab}. This naturally gives rise to the notion of tableaux stabilization. 

\begin{Definition} \label{def:stab} For any skew standard tableau $ T $ with weakly decreasing row sizes from top to bottom, let $ T^{(q)} $ denote the result of attaching $ (q - 1) $ shifted copies of $ T $ to the right of $ T $ so that the result is a skew standard tableau. Let $ n $ denote the size of $ T $ and $ q $ be a positive integer. Let $ \Rect $ denote Sch\"utzenberger's rectification operator \cite{MR0498826}. We say $ T $ \emph{stabilizes} at $ q $ if each element of $[(q - 1)n + 1, qn]$ in $\Rect(T^{(q)})$ lies in the same row as it does in $T^{(q)}$. Let $ \stab(T) $ denote the minimum value at which $ T $ stabilizes.

\end{Definition}

\begin{Remark} Having weakly decreasing row sizes from top to bottom is required for $ T^{(q)} $ to be a skew standard tableau for all $ q \ge 1 $. If this is not the case, $ T^{(q)} $ need not be standard. Consider \Cref{fig:tabincrow}, where $ T^{(2)} $
\begin{figure}[htbp]
    \input{Figures/tabincrow}
\caption{The number of boxes in each row has to be weakly decreasing from top to bottom, or else $ T^{(q)} $ need not be standard.}
\label{fig:tabincrow}
\end{figure}
is not a skew tableau both because of its shape and the third column not being increasing.

\end{Remark}

\begin{Example} \label{ex:stab}
\begin{figure}[htbp]
    \input{Figures/tabexamplepadding1}
\caption{Example of a skew tableau $T$, the tableau $T^{(3)}$ constructed from $T$, and its rectification.}
\label{fig:tabexamplepadding1}
\end{figure}
In Figure~\ref{fig:tabexamplepadding1}, $ 2 $ does not lie in the same row in $ T^{(3)} $ and $ \Rect(T^{(3)}) $, but $ 7, 8, \dots, 12 $ do. Hence, $ \stab(T) = 2 $. As $ 13, 14, \dots, 18 $ also stay in the same row in $ T^{(3)} $ and $ \Rect(T^{(3)}) $, $ T $ stabilizes at 3 as well. In Figure~\ref{fig:tabexamplepadding2}, 9 does not lie in the same row in $ T^{(3)} $ and $ \Rect(T^{(3)}) $, but $ 15, 16, \dots, 21 $ do, so $ \stab(T) = 3 $.

\begin{figure}[htbp]
    \input{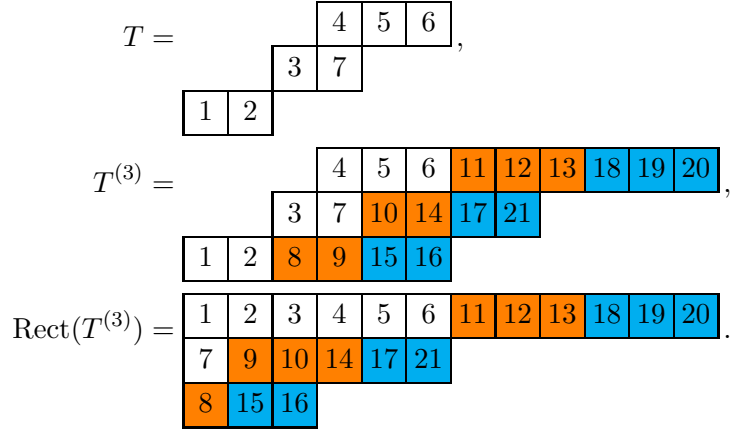}
\caption{Example of a skew tableau $T$, the tableau $T^{(3)}$ constructed from $T$, and its rectification.}
\label{fig:tabexamplepadding2}
\end{figure}

\end{Example}

Ahlbach proved that once a tableaux stabilizes at a value, it continues to stabilize for higher values, \cite[Lemma 3.9]{TabStab}. In the special case where the rows of $ T $ have the same size, he derived a formula for the shape of $ \Rect(T^{(q)}) $ for $ q \ge r - 1 $, \cite[Theorem 1.6]{TabStab}, and used it to show that any tableaux with $ r $ rows of the same size stabilizes at $ r $, \cite[Theorem 1.4]{TabStab}. He conjectured that the same bound still holds when the rows have weakly decreasing sizes but was only able to prove a bound of $ \max(1,2r - 2) $. The main purpose of this paper is to prove this conjecture.

\begin{Theorem} \label{thm:main} For any skew standard tableaux $ T $ with $ r $ rows and weakly decreasing row sizes from top to bottom,
\[
    \stab(T) \le r.
\]
\end{Theorem}

The proof of the same-size-rows case of \Cref{thm:main} in \cite{TabStab}, relied on a formula for the shape of the stabilized tableau. Unfortunately, an analogous formula for the general case would have to involve new terms not present in the previous shape formula, and we have not found such a formula. 

Yet, we generalize the lattice path argument in the proof of Lemma 4.2 in \cite{TabStab} to prove \Cref{thm:main}. Our argument relies on Greene's Theorem \cite{MR0354395} characterizing the shape of the insertion tableaux of words coming from the RSK-correspondence \cite{schensted_1961}, in terms of their increasing subsequences and a careful analysis of the increasing subsequences of reading words of $ T^{(q)} $ for skew standard tableaux $  T $.

In section 2, we show that the family of increasing sequences can be encoded by a family of lattice paths. In section 3, we go over various operations on a family of lattice paths that weakly increase the maximum combined length of corresponding increasing subsequences. In section 4, we prove the main result using the tools developed in the previous sections.

\section{Longest Increasing Subsequences}
In this section we will go over Greene's theorem and set up the language of the matrix and lattice paths we will use. For a more detailed introduction and applications of longest increasing sequences, we recommend the reader to~\cite{romik_2015}. Before we introduce the essential tools, we will briefly explain why we need them.

\begin{figure}[htbp]
    \input{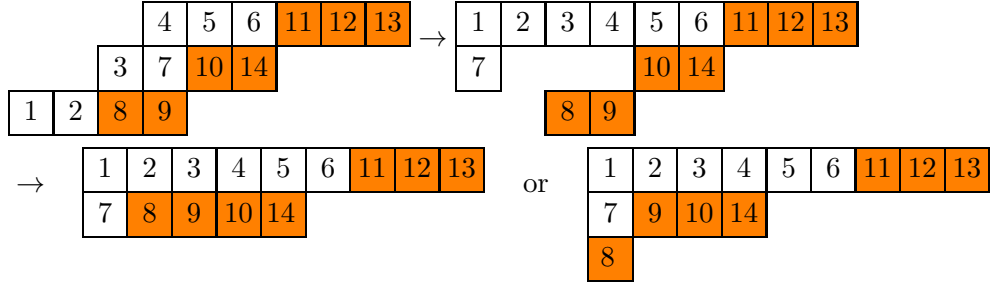}
\caption{If we rectify only the white copy first, rectifying the orange copy can produce two different outcomes depending on the order of the slides.}
\label{fig:obviousbad}
\end{figure}

\begin{Remark}

Rectifying $T^{(q)}$ cannot be done one piece at a time. If one rectifies $T^{(q-1)}$ and then tries to attach another copy of $T$ and apply jeu de taquin, the result need not be unique, unlike rectifying all at once \cite{MR0498826}. Figure~\ref{fig:obviousbad} demonstrates this. Hence we instead look at longest increasing sequences and Greene's theorem to analyze the shape of $\Rect(T^{(q)})$.
\end{Remark}

\subsection{Greene's theorem}

\begin{Definition}(Reading Word)
The \emph{reading word} of a tableau is the word obtained by concatenating the rows from bottom to top. For a non-skew tableau $ T $, let $ \sh(T) $ denote its shape.
\end{Definition}

\begin{figure}[htbp]
    \input{Figures/Tforreading}
\caption{The reading word of $ T $ is $273591468$. The tableau to the right is its rectification, and $ \sh(\Rect(T)) = (5,3,1). $}
\label{fig:Tforreading}
\end{figure}

\begin{Theorem}(Greene's Theorem, \cite{MR0354395})
Let $\pi$ be the reading word of a (skew) standard Young tableau $T$ and let $\ell_k(\pi)$ denote the maximum combined length of $k$ disjoint increasing subsequences of $\pi$. Then, 
\[
    \ell_k(\pi) = \sh(\Rect(T))_1 + \cdots + \sh(\Rect(T))_k.
\]
\end{Theorem}

Consider the tableau in Figure~\ref{fig:Tforreading}. A longest increasing sequence of its reading word $w = 273591468$ is exhibited by $23468$ or $23568$, which fits the fact that the first row has 5 boxes. Two disjoint increasing sequences with the longest combined length are $2359$ and $1468$, which again fits the fact that the combined number of boxes in the first and second row of $ \Rect(T) $ is 8. Note that neither $23468$ nor $23568$ can be reused to obtain two disjoint increasing subsequences of $w $ with total size 8. Finally, $ 27, 359, 1468 $ make 3 disjoint increasing subsequences comprising all of $ w $, so $ \Rect(T) $ has 3 rows.

\begin{Remark}
All the words we use from now on will have no repeated letters.
\end{Remark}


\subsection{Lattice paths}

Suppose $T$ is a skew standard tableau with weakly decreasing row sizes from top to bottom. Let $ n $ denote the size of $ T $, and $r$ the number of rows of $ T $. For each positive integer $j$, let $T + (j - 1)n $ be obtained from $T$ by shifting the entries up by $(j - 1)n $. We define $T^{(q)}$ to be obtained by concatenating $ T, T + n,\ldots,T + (q - 1)n $ together from left to right. Let $ w(q,T) $ denote the reading word of $ T^{(q)} $.

For each positive integer $q$, we create an $r$-by-$q$ matrix $M = M(q,T)$ from $T$: each entry $M_{i,j}$ is a word set as the $(r-j+1)$-th row of $T$, with all entries shifted up by $(i-1)n$. Then, $ w(q,T) $ is the concatenation of $ M(q,T) $ from bottom to top.

\begin{Remark}
We label coordinates Cartesian-style rather than matrix-style since that is more natural for lattice paths. Here, $M_{i,j}$, the $(i,j)$-entry of the matrix $M$, denotes the entry at the $i$-th column from the left and $j$-th row from the bottom. \end{Remark}


\begin{Example} With $ T $ as in \Cref{fig:Tforreading}, we have
 \[
     M(3,T) = \begin{matrix} \color{black}{1468} & \color{orange}{1468} & \color{cyan}{1468} \\ \color{black}{359} & \color{orange}{359} & \color{cyan}{359} \\ \color{black}{27} & \color{orange}{27} & \color{cyan}{27} \end{matrix}
\]
The orange and cyan colors indicate that the entries are shifted up by $9$ and $18$ respectively. So $\color{orange}{3}$ is actually $3+9 = 12$ whereas $\color{cyan}{3}$ is actually $3+18 = 21$. We have $M_{1,1} = 27$ and $M_{3,2} = \color{cyan}{359}$. Note that
\[
    w(3,T) = \color{black}{27} \color{orange}{27} \color{cyan}{27} \color{black}{359} \color{orange}{359} \color{cyan}{359} \color{black}{1468} \color{orange}{1468} \color{cyan}{1468}
\]
is in fact $ M(3,T) $ read bottom to top and concatenated.
\end{Example}

\begin{Definition} For any sequences $ A, B $, let $ A \subseq B $ denote that $ A $ is a subsequence of $ B $. For a sequence $ A $ and a set $I$ we use $A \mid_I$ to denote the restriction of $A$ to the letters in $I$. For example, if $ A = 7164532 $, then
\[
    A \mid_{ \{ 1, 2, 3, 4 \} } = 1432 \subseq A.
\]

\end{Definition}


\begin{Remark}\label{rem:inc}
The words in column $ j $ of $M$ use the letters $I_j \coloneqq \{(j-1)n+1,\ldots,jn\}$ without repetition. Hence any entry used in column $ j $ of $M$ is bigger than any entry used in the previous columns.
\end{Remark}

A \newword{lattice path} within $M$ is a sequence of entries of $M$ that move adjacently right or up at each step. Given a lattice path $S$, we use $S|_{i,*}$ to denote the subpath of $S$ by restricting ourselves to column $i$ of $M$. Similarly we use $S|_{*,j}$ to denote the subpath of $S$ by restricting ourselves to row $ j $ of $M$.

\begin{Definition} \label{def:mclis} A collection of words is disjoint if no two words have any letters in common. For any words $ A_1, \dots, A_k $, let
\begin{align*}
    \ell(A_1, \dots, A_k) = & \tx{ the maximum combined length of disjoint} \\
            & \tx{ increasing subsequences of } A_1, \dots, A_k \tx{ respectively}.
\end{align*}
We say that a collection of disjoint increasing sequences $S_1 \subseq A_1,\ldots,S_k \subseq A_k$ \newword{exhibits} $\ell(A_1,\ldots,A_k)$ if $|S_1|+\cdots+|S_k|=\ell(A_1,\ldots,A_k)$.
\end{Definition}


For example, $ \ell(152643,13497) = 6 $ is exhibited by the pairs $(26, 1347) $ or $( 156, 349 )$. Any letter, like  $ 1 $ in particular, cannot be used by both subsequences since the subsequences have to be disjoint.

\bs
We will often identify lattice paths with the concatenated sequence of entries in $ M $ they represent as we will do in \Cref{lem:LPsubseq} and \Cref{fig:latticepath1}. We are investigating these lattice paths because we are interested in the increasing subsequences they contain.

\begin{Lemma} \label{lem:LPsubseq}
Any increasing subsequence of $ w(q,T) $ is a subsequence of a lattice path within $M$.
\end{Lemma}

\begin{proof}
Each increasing subsequence of $ w(q,T) $ cannot go downward in $ M $ because that is going backward in reading order of $ T^{(q)} $. It also cannot go left in $ M $ because any entry in column $ j $ of $M$ is bigger than any entry in column $ (j - 1) $ of $ M $, by \Cref{rem:inc}. Finally, each entry of $ M $ is an increasing word, so increasing subsequences can only travel right within these points as well. Therefore any increasing subsequence of $ w(q,T) $ only travels right and up in $ M $ and thus must be a subsequence of a lattice path. 

\end{proof}

By \Cref{lem:LPsubseq}, we can view any $ k $ disjoint increasing subsequences of $ w(q,T) $ as contained within a family of $ k $ possibly intersecting lattice paths. And for a family $ P_1, \dots, P_k $ of lattice paths, we let $ \ell(P_1, \dots, P_k) $ denote the maximum combined length of disjoint increasing subsequences of $ P_1, \dots, P_k $, respectively, as in \Cref{def:mclis}.

\begin{Example}

With $ T $ as in \Cref{fig:Tforreading}, we have
\[
    \Rect(T^{(3)}) = \byt 1 & 3 & 4 & 6 & 8 & *(orange) 10 & *(orange) 13  & *(orange) 15 & *(orange) 17 & *(cyan) 19 & *(cyan) 22 & *(cyan) 24 & *(cyan) 26 \\ 2 & 5 & 9 & *(orange) 12 & *(orange) 14 & *(orange) 18 & *(cyan) 21 & *(cyan) 23 & *(cyan) 27 \\ 7 & *(orange) 11 & *(orange) 16 & *(cyan) 20 & *(cyan) 25 \eyt 
\]
with shape $ \lam = (13, 9, 5) $ and
\[
     M(3,T) = \begin{matrix} \color{black}{1468} & \color{orange}{1468} & \color{cyan}{1468} \\ \color{black}{359} & \color{orange}{359} & \color{cyan}{359} \\ \color{black}{27} & \color{orange}{27} & \color{cyan}{27} \end{matrix} \, .
\]
A longest increasing subsequence of $ w(3,T) $ is $ 23568\color{orange}{1468} \color{cyan}{1468} $, which has size 13 and is a subsequence of the lattice path in \Cref{fig:latticepath1}.
\begin{figure}[htbp]
    \[	
	\begin{tikzpicture}
\draw (0,0) -- (0,2) -- (2,2);
\draw [fill = black] (0,0) circle (0.050);
\draw [fill = black] (0,1) circle (0.050);
\draw [fill = black] (0,2) circle (0.050);
\draw [fill = black] (1,0) circle (0.050);
\draw [fill = black] (1,1) circle (0.050);
\draw [fill = black] (1,2) circle (0.050);
\draw [fill = black] (2,0) circle (0.050);
\draw [fill = black] (2,1) circle (0.050);
\draw [fill = black] (2,2) circle (0.050);
\end{tikzpicture}
\]
\caption{The lattice path 273591468\color{orange}{1468}\color{cyan}{1468} \color{black} containing $ 23568\color{orange}{1468} \color{cyan}{1468} $.}
\label{fig:latticepath1}
\end{figure}
Note this agrees with Greene's Theorem since $ \lam_1 = 13 $. Moreover, two disjoint increasing subsequences of $ w(3,T) $ with the longest combined length are $ 1468\color{orange}{1468} \color{cyan}{1468} $ and $ 2359\color{orange}{359} \color{cyan}{359} $, have total size 22 and are subsequences of the lattice paths in \Cref{fig:latticepath2}.
\begin{figure}[htbp]
    \[	
	\begin{tikzpicture}
\draw (0,2) -- (2,2);
\draw (0,0) -- (0,1) -- (2,1);
\draw [fill = black] (0,0) circle (0.050);
\draw [fill = black] (0,1) circle (0.050);
\draw [fill = black] (0,2) circle (0.050);
\draw [fill = black] (1,0) circle (0.050);
\draw [fill = black] (1,1) circle (0.050);
\draw [fill = black] (1,2) circle (0.050);
\draw [fill = black] (2,0) circle (0.050);
\draw [fill = black] (2,1) circle (0.050);
\draw [fill = black] (2,2) circle (0.050);
\end{tikzpicture}
\]
\caption{The minimal lattice paths containing $ 1468\color{orange}{1468} \color{cyan}{1468} $ and $ 2359\color{orange}{359} \color{cyan}{359} $.}
\label{fig:latticepath2}
\end{figure}
Again, this agrees with Greene's Theorem since $ \lam_1 + \lam_2 = 13 + 9 = 22 $. Finally, three disjoint increasing subsequences of $ w(3,T) $ with the longest combined length are $ 1468\color{orange}{1468} \color{cyan}{1468}, \color{black}{359}\color{orange}{359} \color{cyan}{359}, $ and $ 27\color{orange}{27} \color{cyan}{27} $, which have total size 27 and are subsequences of the lattice paths in \Cref{fig:latticepath3}.
\begin{figure}[htbp]
    \[	
	\begin{tikzpicture}
\draw (0,2) -- (2,2);
\draw (0,1) -- (2,1);
\draw (0,0) -- (2,0);
\draw [fill = black] (0,0) circle (0.050);
\draw [fill = black] (0,1) circle (0.050);
\draw [fill = black] (0,2) circle (0.050);
\draw [fill = black] (1,0) circle (0.050);
\draw [fill = black] (1,1) circle (0.050);
\draw [fill = black] (1,2) circle (0.050);
\draw [fill = black] (2,0) circle (0.050);
\draw [fill = black] (2,1) circle (0.050);
\draw [fill = black] (2,2) circle (0.050);
\end{tikzpicture}
\]
\caption{The minimal lattice paths containing $ 1468\color{orange}{1468} \color{cyan}{1468}, \color{black}{359}\color{orange}{359} \color{cyan}{359}, $ and $ 27\color{orange}{27} \color{cyan}{27} $.}
\label{fig:latticepath3}
\end{figure}
Again, this agrees with Greene's Theorem since $ \lam_1 + \lam_2 + \lam_3 = 13 + 9 + 5 = 27 $.


\end{Example}
\color{black}

\subsection{Splitting words} In this subsection, we show that we can attack each column of a family of lattice paths separately. We also prove a property of shuffling two increasing sequences that will be a key tool in the following section.

\begin{Lemma}\label{lem:ColumnSplitting}(Column Splitting) For any lattice paths $P_1,\dots,P_k$, we have
\[
    \ell(P_1,\ldots,P_k) = \sum_{j = 1}^q \ell(P_1 \mid_{j,*},\ldots,P_k \mid_{j,*}).
\]
\end{Lemma} 

\begin{proof}
Since any element of column $ j $ is larger than any element of column ${j-1}$ for all $ j $ (\Cref{rem:inc}), the increasing subsequences of each path can be picked independently column by column. 
\end{proof}

Hence, when analyzing a family of lattice paths and its length $\ell(S_1,\ldots,S_m)$, we may study each column separately.

\begin{Definition} \label{def:splitpt} A \textbf{splitting point} of a word $ w = w_1 \dots w_n $ is an index $ k $ where 
\[
    w_i < w_j \quad \text{ for all } i \leq k < j.
\]
\end{Definition}

\begin{Lemma}\label{lem:splpoint} Suppose $ w $ is a word that is a shuffle of two increasing sequences, where one of these subsequences includes both the initial and final letters. Then $ w $ has a splitting point.
\end{Lemma}

\begin{proof}
Without loss of generality, we may standardize $w = w_1 \dots w_n $ so it uses $ [n] $. By assumption, $ w $ is a shuffle of two increasing subsequences, $ A $ and $ B $, where $ A $ includes $ w_1 $ and $ w_n $. We wish to show that $w$ has a splitting point. 

If $w_n=n$, then $n-1$ is a splitting point. Thus we may assume that $w_n < n$, which implies $ \max(A) = w_n < n $, so $ n \in B$.

For an index $i$, we say $ i $ is $A$\emph{-winning} if $\max(w_1, w_2, \dots, w_i) \in A$ and $B$\emph{-winning} otherwise. We have that $1$ is $A$-winning because $ w_1 \in A $, and $n$ is $B$-winning since $n \in B$.

Thus, there exists a  smallest positive integer $ k $ such that $k$ is B-winning. Since $k \geq 2$, $k-1$ must be A\emph{-winning}. It follows that $\max(w_1, w_2, \dots, w_{k-1}) \in A $ and $ w_k = \max(w_1, w_2 \dots, w_k) \in B $.

We claim that $k-1$ is a splitting point of $w$. Let $i$ and $j$ be positive integers such that $i \leq k-1 < j$. If $ w_j \in A $, then 
\[
    w_i \le \max(w_1, w_2, \dots, w_{k-1}) < w_j
\]
since $\max(w_1, w_2, \dots, w_{k-1}) \in A $ and $ A $ is increasing. If $ w_j \in B $, then
\[
    w_i < \max(w_1, w_2 \dots, w_k) = w_k \leq w_j
\]
since $ w_k \in B $ and $ B $ is increasing. Thus, by \Cref{def:splitpt}, $k-1$ is a splitting point of $ w $.
\end{proof}

For example, the word $ w = 231547968 $, which is a shuffle of $ 23568 $ and $ 1479 $, has a splitting point at $ 5 $ since $ 2, 3, 1, 5, 4 < 7, 9, 6, 8 $.


\section{Lattice Path Tools} 
\label{sec:LPtools}
In this section we will introduce and prove various tools we will need for the proof of our main result. Given a family of paths, we will show several operations that modify those paths while weakly increasing the maximum combined length of increasing subsequences they contain. 

  The first tool, top-down switching, will allow us to modify the paths so that they do not cross (intersections can happen, but they will be non-transversal), and thus can be ordered in a top-down fashion. The second tool, left-shifting, will allow us to shift portions of vertical segments of a path to the left as long as no new intersections appear.  The last two tools, rectangular and reverse rectangular flip, will allow us to split paths with shared horizontal segments, while maintaining a top-down order.




\subsection{Top-Down Switching}

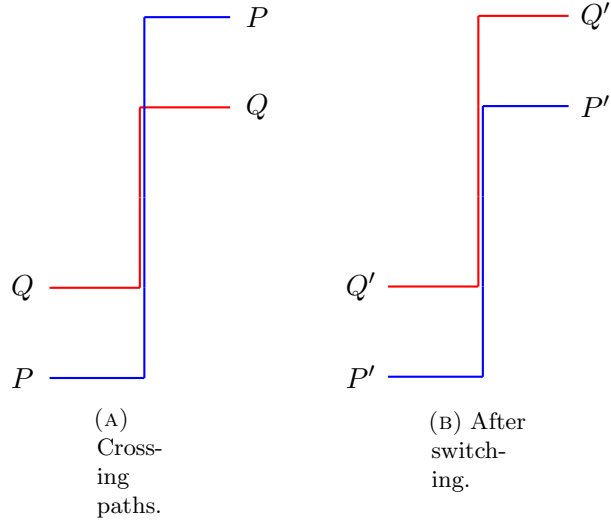
\begin{figure}
        \centering
        \subfloat[Crossing paths.]{\begin{tikzpicture}[scale = 1.2]
\coordinate (-1) at (0,0) {};
\coordinate (0) at (0,1) {};
\coordinate (1) at (1, 1) {};
\coordinate (1') at (1,2) {};
\coordinate (2) at (1,3) {};
\coordinate (3) at ( 2,3 ) {};
\coordinate (4) at ( 1.05,0 ) {};
\coordinate (5) at ( 1.05,1 ) {};
\coordinate (6) at ( 1.05,2 ) {};
\coordinate (7) at ( 1.05,3 ) {};
\coordinate (8) at ( 1.05,4 ) {};
\coordinate (9) at (2,4){};

\node at (-0.3,0) {$ P $};
\node at (-0.3,1) {$ Q $};

\node at (2.3,3) {$ Q $};
\node at (2.3,4) {$ P $};

\draw [red,thick](0) -- (1) node   {};
\draw [red,thick](2) -- (3) node   {};
\draw [red,thick](1) -- (1') node   {};
\draw [red,thick](1') -- (2) node   {};
\draw [blue,thick](-1) -- (4) node {};
\draw [blue,thick](4) -- (5) node  {};
\draw [blue,thick](5) -- (6) node {} ;
\draw [blue,thick](6) -- (7) node  {};
\draw [blue,thick](7) -- (8) node {} ;
\draw [blue,thick](8) -- (9) node {};
\end{tikzpicture}}
        \qquad
        \subfloat[After switching.]{\begin{tikzpicture}[scale = 1.2]

\coordinate (-1) at (0,0) {};
\coordinate (0) at (0,1) {};
\coordinate (1) at (1, 1) {};
\coordinate (1') at (1,2) {};
\coordinate (2) at (1,3) {};
\coordinate (3) at ( 2,3 ) {};
\coordinate (4) at ( 1.05,0 ) {};
\coordinate (5) at ( 1.05,1 ) {};
\coordinate (6) at ( 1.05,2 ) {};
\coordinate (7) at ( 1.05,3 ) {};
\coordinate (8) at ( 1,3 ) {};
\coordinate (9) at ( 1,4 ) {};
\coordinate (10) at ( 2,4 ) {};

\node at (-0.3,0) {$ P' $};
\node at (-0.3,1) {$ Q' $};

\node at (2.3,3) {$ P' $};
\node at (2.3,4) {$ Q' $};

\draw [red,thick](0) -- (1) node   {};
\draw [blue,thick](7) -- (3) node   {};
\draw [red,thick](1) -- (1') node   {};
\draw [red,thick](1') -- (2) node   {};
\draw [red,thick](9) -- (10) node   {};

\draw [blue,thick](-1) -- (4) node  {};
\draw [blue,thick](4) -- (5) node  {};
\draw [blue,thick](5) -- (6) node {} ;
\draw [blue,thick](6) -- (7) node  {};
\draw [red,thick](8) -- (9) node {} ;

\end{tikzpicture}}
    \caption{The top-down switching process.}
    \label{fig:topdown}
    \end{figure}

Consider two lattice paths $ P $(blue) and $ Q $(red) within the matrix $ M $ that cross in a way that they switch from $ Q $ being on top to $ P $ being on top, as in the left figure of Figure~\ref{fig:topdown}.

Consider the result of switching the labels on the paths after $ P $ and $ Q $ diverge, giving new paths $ P' $ (blue) and $ Q' $ (red) as shown in the right figure of Figure~\ref{fig:topdown}. After this top-down switch, we have a clear upper and lower path. Our goal in this subsection is to show that this operation weakly increases the maximum combined length of increasing subsequences they contain.

\begin{Lemma} \label{lem:LIS1} For any words $ A, B, C, D_1, \dots, D_m $, 
\[
    \ell(ABC,B,D_1, \dots, D_m) \leq \ell(AB,BC,D_1, \dots, D_m).
\]
\end{Lemma}

\begin{proof}

Let $S_{ABC},S_{B},\ldots$ be a list of sequences that exhibits $l(ABC,B,D_1,\ldots,D_m)$. We will modify $S_{ABC}$ and $S_{B}$ while preserving their union, into disjoint increasing subsequences of $AB$ and $BC$ respectively, making the remaining disjoint increasing subsequences irrelevant. To do this, consider $S = S_{ABC} \sqcup S_B $ with its elements ordered how they appear in $ ABC $ so that $ S \subseq ABC $.

If $S_{ABC}$ has no elements from $A$ or no elements from $C$, we can assign $ S_{ABC} $ as a subsequence of $ BC $ or $ AB $, respectively and then $ S_B $ as a subsequence of the other. Otherwise, $ S _{ABC} $ has an element from $ A $ and an element from $C $, so $ S_{ABC} $ contains the first and last elements of $ S $. As $ S $ is a shuffle of $ S_{ABC} $ and $ S_B $, $ S $ has a splitting point by Lemma~\ref{lem:splpoint}. If $ S_{ABC} $ contains the maximum element of $ B $ in $ S $, then we can put the splitting point right after $ B $. Else, by Lemma~\ref{lem:splpoint}, there exists a splitting point where the sequence first switches from $ S_{ABC} $-winning to $ S_{B} $-winning, which is right before some element of $ B $.  

Next, we can partition $S_{ABC}$ as
\[
    S_{ABC} = S_{ABC}^{1} S_{ABC}^2
\]
where $ S_{ABC}^{1} $ is the subsequence of $ S_{ABC} $ before the splitting point, and $ S_{ABC}^{2} $ is the subsequence after the splitting point. Similarly, split $ S_B $ into 
\[
    S_B = S_B^1 S_B^2
\]
by splitting at the splitting point. Now create the sequences 
\begin{align*}
    S_{AB} = S_{ABC}^{1} S_B^2 \subseq AB, \\
    S_{BC} = S_{B}^{1} S_{ABC}^2 \subseq BC.
\end{align*}
These are subsequences of $ AB, BC $ respectively because our splitting point was either at the end of $ B $ or just before some element of $ B $. Furthermore, $ S_{AB}$ and $S_{BC}$ are increasing because every element left of the the splitting point is less than every element right of the splitting point.

Hence we have successfully modified $S_{ABC}$ and $S_{B}$ into $S_{AB} \subseq AB $ and $S_{BC} \subseq BC$ while maintaining their union, which proves \Cref{lem:LIS1}.
\end{proof}
\begin{Example}
Let $ A = 13 , B = 52947, C = 68 $. Then, $ \ell(ABC,B) = 8 $ is exhibited by
\[
    S_{ABC} = 13568, \qquad S_B = 247.
\]
Viewing their union as a subsequence of $ ABC $ gives $ S = 13524768 $. Then, since $ S $ has a splitting point at 5, just after the 4, we split them into
\begin{align*}
    S_{ABC}^1 & = 135, \qquad S_{ABC}^{2} = 68, \\
    S_B^1 & = 24, \qquad S_B^2 = 7 
\end{align*}
and then recombine them into
\begin{align*}
    S_{AB} = S_{ABC}^1 S_B^2 = 1357 \subseq AB, \\
    S_{BC} = S_B^1 S_{ABC}^2 = 2468 \subseq BC.
\end{align*}
This shows $ \ell(AB,BC) \ge 8 = \ell(ABC,B) $.
\end{Example}

\begin{Lemma}(Top-Down Switching) \label{lem:TDS}
Suppose $ P, Q $ are lattice paths, and there exists a column $ j $ where $ P $ begins weakly lower in column $ j $ than $ Q $ but $ P $ ends weakly higher in column $ j $ than $ Q $. Then, performing the top-down swap as in \Cref{fig:topdown},
\[
    \ell(P_1, \dots, P_k, P,Q) \le \ell(P_1, \dots, P_k, P',Q'). 
\]
for any lattice paths $ P_1, \dots,P_k $.
\end{Lemma}

\begin{proof}

By the Column Splitting lemma (Lemma~\ref{lem:ColumnSplitting}), we only need to compare the lengths in the column $ j $ where $ P $ and $ Q $ intersect since in every other column, we are measuring the maximum combined length of disjoint increasing subsequences of exactly the same subsequences. Using the labels in
\Cref{fig:topdownlabels}: 
\begin{enumerate}
    \item $ A $ is the sequence in column $ j $ strictly below the intersection,
    \item $ B $ is the sequence in column $ j $ contained in the intersection of $ P $ and $ Q $,
    \item $ C $ is the sequence in column $ j $ strictly above the intersection,
    \item $ D_1, \dots, D_m $ are the sequences in the non-intersecting columns that are common to both families of paths,
\end{enumerate}

\begin{figure}
        \centering
        \subfloat[]{\begin{tikzpicture}[scale = 1.2]
\coordinate (-1) at (0,0) {};
\coordinate (0) at (0,1) {};
\coordinate (1) at (1, 1) {};
\coordinate (1') at (1,2) {};
\coordinate (2) at (1,3) {};
\coordinate (3) at ( 2,3 ) {};
\coordinate (4) at ( 1.05,0 ) {};
\coordinate (5) at ( 1.05,1 ) {};
\coordinate (6) at ( 1.05,2 ) {};
\coordinate (7) at ( 1.05,3 ) {};
\coordinate (8) at ( 1.05,4 ) {};
\coordinate (9) at (2,4){};

\node at (-0.3,0) {$ P $};
\node at (-0.3,1) {$ Q $};

\node at (0.8,0.5) {$ A $};
\node at (0.8,2) {$ B $};
\node at (0.8,3.5) {$ C $};

\node at (2.3,3) {$ Q $};
\node at (2.3,4) {$ P $};

\draw [red,thick](0) -- (1) node   {};
\draw [red,thick](2) -- (3) node   {};
\draw [red,thick](1) -- (1') node   {};
\draw [red,thick](1') -- (2) node   {};
\draw [blue,thick](-1) -- (4) node {};
\draw [blue,thick](4) -- (5) node  {};
\draw [blue,thick](5) -- (6) node {} ;
\draw [blue,thick](6) -- (7) node  {};
\draw [blue,thick](7) -- (8) node {} ;
\draw [blue,thick](8) -- (9) node {};
\end{tikzpicture}}
        \qquad
        \subfloat[]{\begin{tikzpicture}[scale = 1.2]

\coordinate (-1) at (0,0) {};
\coordinate (0) at (0,1) {};
\coordinate (1) at (1, 1) {};
\coordinate (1') at (1,2) {};
\coordinate (2) at (1,3) {};
\coordinate (3) at ( 2,3 ) {};
\coordinate (4) at ( 1.05,0 ) {};
\coordinate (5) at ( 1.05,1 ) {};
\coordinate (6) at ( 1.05,2 ) {};
\coordinate (7) at ( 1.05,3 ) {};
\coordinate (8) at ( 1,3 ) {};
\coordinate (9) at ( 1,4 ) {};
\coordinate (10) at ( 2,4 ) {};

\node at (-0.3,0) {$ P' $};
\node at (-0.3,1) {$ Q' $};

\node at (0.8,0.5) {$ A $};
\node at (0.8,2) {$ B $};
\node at (0.8,3.5) {$ C $};

\node at (2.3,3) {$ P' $};
\node at (2.3,4) {$ Q' $};

\draw [red,thick](0) -- (1) node   {};
\draw [blue,thick](7) -- (3) node   {};
\draw [red,thick](1) -- (1') node   {};
\draw [red,thick](1') -- (2) node   {};
\draw [red,thick](9) -- (10) node   {};

\draw [blue,thick](-1) -- (4) node  {};
\draw [blue,thick](4) -- (5) node  {};
\draw [blue,thick](5) -- (6) node {} ;
\draw [blue,thick](6) -- (7) node  {};
\draw [red,thick](8) -- (9) node {} ;

\end{tikzpicture}}
    \caption{Labels on the top-down switching paths}
    \label{fig:topdownlabels}
    \end{figure}
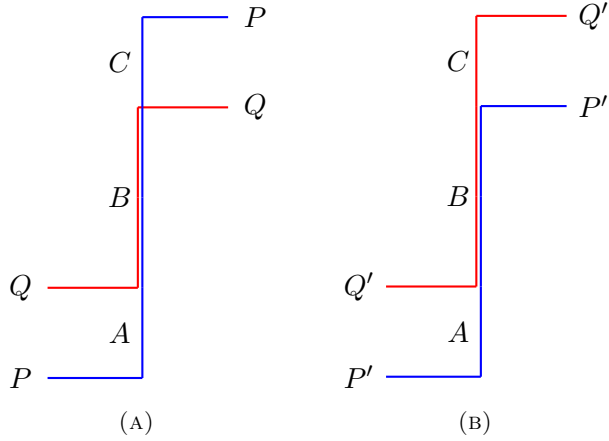
    
we have
\begin{align*}  
    \ell(P_1, \dots, P_k, P,Q) & = \ell(ABC,B,D_1, \dots, D_m)
\end{align*}  
before the top-down switch and
\begin{align*}
     \ell(P_1, \dots, P_k, P',Q') = \ell(AB,BC,D_1, \dots, D_m)
\end{align*}
after the top-down switch. The fact that the remaining portions of $ P $ and $ Q $ after column $ j $ switch in $ P' $ and $ Q' $ is irrelevant by Column Splitting. Thus,
\begin{align*}
    \ell(P_1, \dots, P_k, P,Q) & = \ell(ABC,B,D_1, \dots, D_m) \\ 
        & \leq \ell(AB,BC,D_1, \dots, D_m) \qquad \tx{ by \Cref{lem:LIS1} }  \\
        & = \ell(P_1, \dots, P_k, P',Q').
\end{align*}

\end{proof}

\begin{figure}
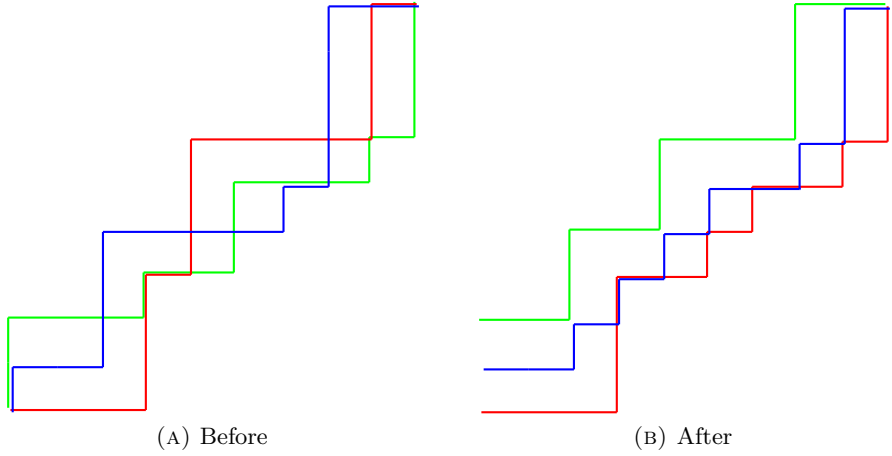

        \centering
        \subfloat[Before]{\input{Figures/newrectflip}}
        \qquad
        \subfloat[After]{\input{Figures/newtopdown}}
        
    \caption{After applying Top-Down Switching, we can label the newly ordered paths $P_1,P_2,P_3$ from top to bottom, in a way that a path $P_i$ always stays weakly above $P_{i+1}$.}
    \label{fig:newtopdown}
    \end{figure}

An example of top-down switching is given in Figure~\ref{fig:newtopdown}. Each application of top-down switching eliminates a transverse crossing from our set of paths. Thus, after applying top-down switching as many times as possible in any order, there are no longer any transverse crossings in the resulting family of paths. Then we can order the family of paths from top to bottom. Hence, top-down switching lets us adjust the paths so they come in a top-down order $ P_1, \dots, P_k $ where $ P_i $ lies weakly above $ P_{i + 1} $ for all $ i $.

\subsection{Left-shifting}

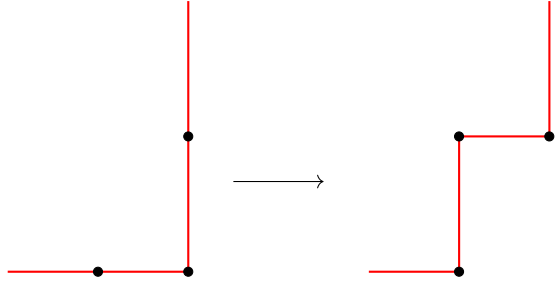
\begin{figure}[htbp]

\begin{center}
    \begin{tikzpicture}[scale = 1.2]
\coordinate (-1) at (-1,0) {};
\coordinate (0) at (0,0) {};
\coordinate (1) at (1, 0) {};
\coordinate (2) at (1,3) {};

\coordinate (-1') at (3,0) {};
\coordinate (0') at (4,0) {};
\coordinate (1') at (4, 1.5) {};
\coordinate (2') at (5, 1.5) {};
\coordinate (3') at (5, 3) {};

\draw [red,thick](-1) -- (0) node   {};
\draw [red,thick](0) -- (1) node   {};
\draw [red,thick](1) -- (2) node   {};

\draw [red,thick](-1') -- (0') node {};
\draw [red,thick](0') -- (1') node {};
\draw [red,thick](1') -- (2') node  {};
\draw [red,thick](2') -- (3') node  {};

\draw [->] (1.5,1) -- (2.5,1);

\draw [fill = black] (0,0) circle (0.050);
\draw [fill = black] (1,0) circle (0.050);
\draw [fill = black] (4,0) circle (0.050);
\draw [fill = black] (1,1.5) circle (0.050);
\draw [fill = black] (4,1.5) circle (0.050);
\draw [fill = black] (5,1.5) circle (0.050);

\end{tikzpicture}
\end{center}
\caption{Left-shifting applied to a lattice path}
\label{fig:leftshift1}
\end{figure}

\begin{Definition}
A \newword{left-shift} of a lattice path is an operation where part of the column of a path, following a horizontal segment of length at least two, is shifted one column to the left and creates no new intersections, as in Figure \ref{fig:leftshift1}.
\end{Definition}

We are introducing left-shifting because it again weakly increases the maximum combined length of increasing subsequences contained in the lattice paths.

\begin{Lemma} \label{lem:leftshift} Suppose $ P_1, \dots, P_k, P $ are a family of lattice paths where a left-shift can be applied to $ P $ in a way so that no new intersection is introduced. Letting $ Q $ denote this left-shift of $ P $, we have
\[
    \ell(P_1, \dots, P_k, P) \le \ell(P_1, \dots, P_k, Q).
\]
\end{Lemma}

\begin{proof} Suppose $ \ell(P_1, \dots, P_k, P) $ is exhibited by $ S_1, \dots, S_k, S $. Now, we restrict our attention only to the column of $ P $ that was just shifted, say $ A $ is sequence of the bottom point in the shifted column, $ B $ is the sequence strictly between the boundary points of the shifted piece, $ C $ is the sequence at the top point of the shift, and $ D $ is the sequence strictly above the shift. Also, let $ A', B', C' $ denote $ A, B , C $ decremented by $ n $ component-wise, respectively. Because the row sizes of $ T $ weakly decrease from top to bottom, $ A, C $ are increasing sequences with $ \ell(C) \ge \ell(A) $. See \Cref{fig:leftshift2}.

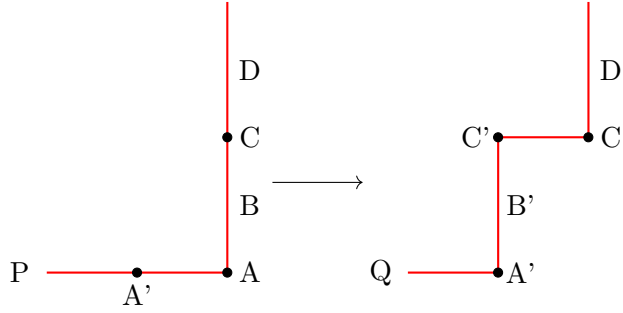
\begin{figure}[htbp]

\begin{center}
    \begin{tikzpicture}[scale = 1.2]
\coordinate (-1) at (-1,0) {};
\coordinate (0) at (0,0) {};
\coordinate (1) at (1, 0) {};
\coordinate (2) at (1,3) {};

\coordinate (-1') at (3,0) {};
\coordinate (0') at (4,0) {};
\coordinate (1') at (4, 1.5) {};
\coordinate (2') at (5, 1.5) {};
\coordinate (3') at (5, 3) {};

\draw [red,thick](-1) -- (0) node   {};
\draw [red,thick](0) -- (1) node   {};
\draw [red,thick](1) -- (2) node   {};

\draw [red,thick](-1') -- (0') node {};
\draw [red,thick](0') -- (1') node {};
\draw [red,thick](1') -- (2') node  {};
\draw [red,thick](2') -- (3') node  {};

\draw [->] (1.5,1) -- (2.5,1);

\node at (1.25, 2.25) {D};
\node at (5.25, 2.25) {D};
\node at (1.25, 1.5) {C};
\node at (5.25, 1.5) {C};
\node at (3.75, 1.5) {C'};
\node at (1.25, 0.75) {B};
\node at (4.25, 0.75) {B'};
\node at (1.25, 0) {A};
\node at (4.25, 0) {A'};
\node at (0, -0.25) {A'};
\node at (-1.3, 0) {P};
\node at (2.7, 0) {Q};

\draw [fill = black] (0,0) circle (0.050);
\draw [fill = black] (1,0) circle (0.050);
\draw [fill = black] (4,0) circle (0.050);
\draw [fill = black] (1,1.5) circle (0.050);
\draw [fill = black] (4,1.5) circle (0.050);
\draw [fill = black] (5,1.5) circle (0.050);

\end{tikzpicture}
\end{center}
\caption{Labeled sections on the lattice path}
\label{fig:leftshift2}
\end{figure}

Let $ S \mid_{AB} = a_1 \dots a_r b_1 \dots b_s $ and $ S \mid_C = c_1 \dots c_t $. Without loss of generality, we may assume $ S $ contains $ A' $. If not, just add the rest of the elements of $ A' $, which preserves $ S $ being increasing as $ A' $ is alone its own column within $ P $. Then, to get an increasing subsequence of $ Q $ of weakly longer length, first remove $ A' $. Next, shift the letters of $ S \mid_{AB} $ to letters of $ A'B' $ by subtracting $n$. For notational convenience, we let $ x' = x - n $. Then, between $ S \mid_{AB} - n $ and $ S \mid_{CD} $, insert $ x' (x' + 1) ... \max(C') \min(C) (\min(C) + 1)... (x - 1) $ where $ x $ is the minimum element of $ C $ greater than $ b_s $ if it exists. Note that by construction $ b_s < x \le c_1 $ because $ b_s < c_1 $. If no such $ x $ exists, which is when $ b_s > \max(C) $, insert $ C $ instead. Either way, we are inserting an increasing sequence of size $ \ell(C) $. Thus,
\[
    S = \dots A' \; a_1 \dots a_r b_1 ... b_s c_1 c_2 ... c_t y \dots
\]
becomes
\begin{align*}
    R = \begin{cases} \dots a_1' \dots a_r' b_1' \dots b_s' x' ... \max(C') \min(C)... (x - 1) c_1 \dots c_t y \dots & \tx{ if } b_s < \max(C) \\
    \dots a_1' \dots a_r' b_1' \dots b_s' C y \dots & \tx{ if } b_s > \max(C).
    \end{cases}
\end{align*}
If $ b_s < \max(C) $, this sequence is increasing since $ b_s < x \le c_1 $. If $ b_s > \max(C) $, the sequence is still increasing since $ \max(C) < b_s < y $. Thus, we have found an increasing $ R \subseq Q $ with length
\[
    \ell(R) = \ell(S) - \ell(A) + \ell(C) \ge \ell(S) \qquad \tx{ because } \ell(C) \ge \ell(A).
\]
Because the space we shifted $ P $ into was empty earlier, $ R $ is still disjoint from $ S_1, \dots, S_k $. Therefore, 
\[
    \ell(P_1, \dots, P_k, P) \le \ell(P_1, \dots, P_k, Q).
\]
\end{proof}

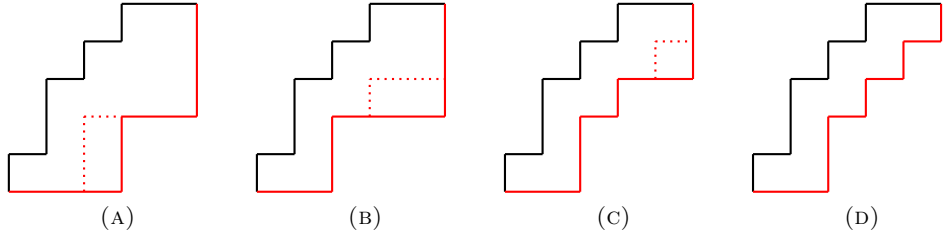
\begin{figure}[htbp]
        \centering
        \subfloat[]{\begin{tikzpicture}[scale = 0.5]
        \coordinate (1a) at (1,1) {};
    \coordinate (1b) at (1,2) {};
    \coordinate (1c) at (1,3) {};
    \coordinate (1d) at (1,4) {};
    \coordinate (1e) at (1,5) {};
    \coordinate (1f) at (1,6) {};
    \coordinate (1g) at (1,7) {};
    \coordinate (1h) at (1,8) {};
    \coordinate (1i) at (1,9) {};
    \coordinate (1j) at (1,10) {};
    \coordinate (2a) at (2,1) {};
    \coordinate (2b) at (2,2) {};
    \coordinate (2c) at (2,3) {};
    \coordinate (2d) at (2,4) {};
    \coordinate (2e) at (2,5) {};
    \coordinate (2f) at (2,6) {};
    \coordinate (2g) at (2,7) {};
    \coordinate (2h) at (2,8) {};
    \coordinate (2i) at (2,9) {};
    \coordinate (2j) at (2,10) {};
    \coordinate (3a) at (3,1) {};
    \coordinate (3b) at (3,2) {};
    \coordinate (3c) at (3,3) {};
    \coordinate (3d) at (3,4) {};
    \coordinate (3e) at (3,5) {};
    \coordinate (3f) at (3,6) {};
    \coordinate (3g) at (3,7) {};
    \coordinate (3h) at (3,8) {};
    \coordinate (3i) at (3,9) {};
    \coordinate (3j) at (3,10) {};
    \coordinate (4a) at (4,1) {};
    \coordinate (4b) at (4,2) {};
    \coordinate (4c) at (4,3) {};
    \coordinate (4d) at (4,4) {};
    \coordinate (4e) at (4,5) {};
    \coordinate (4f) at (4,6) {};
    \coordinate (4g) at (4,7) {};
    \coordinate (4h) at (4,8) {};
    \coordinate (4i) at (4,9) {};
    \coordinate (4j) at (4,10) {};
    \coordinate (5a) at (5,1) {};
    \coordinate (5b) at (5,2) {};
    \coordinate (5c) at (5,3) {};
    \coordinate (5d) at (5,4) {};
    \coordinate (5e) at (5,5) {};
    \coordinate (5f) at (5,6) {};
    \coordinate (5g) at (5,7) {};
    \coordinate (5h) at (5,8) {};
    \coordinate (5i) at (5,9) {};
    \coordinate (5j) at (5,10) {};
    \coordinate (6a) at (6,1) {};
    \coordinate (6b) at (6,2) {};
    \coordinate (6c) at (6,3) {};
    \coordinate (6d) at (6,4) {};
    \coordinate (6e) at (6,5) {};
    \coordinate (6f) at (6,6) {};
    \coordinate (6g) at (6,7) {};
    \coordinate (6h) at (6,8) {};
    \coordinate (6i) at (6,9) {};
    \coordinate (6j) at (6,10) {};
    \coordinate (7a) at (7,1) {};
    \coordinate (7b) at (7,2) {};
    \coordinate (7c) at (7,3) {};
    \coordinate (7d) at (7,4) {};
    \coordinate (7e) at (7,5) {};
    \coordinate (7f) at (7,6) {};
    \coordinate (7g) at (7,7) {};
    \coordinate (7h) at (7,8) {};
    \coordinate (7i) at (7,9) {};
    \coordinate (7j) at (7,10) {};
    \coordinate (8a) at (8,1) {};
    \coordinate (8b) at (8,2) {};
    \coordinate (8c) at (8,3) {};
    \coordinate (8d) at (8,4) {};
    \coordinate (8e) at (8,5) {};
    \coordinate (8f) at (8,6) {};
    \coordinate (8g) at (8,7) {};
    \coordinate (8h) at (8,8) {};
    \coordinate (8i) at (8,9) {};
    \coordinate (8j) at (8,10) {};
    \coordinate (9a) at (9,1) {};
    \coordinate (9b) at (9,2) {};
    \coordinate (9c) at (9,3) {};
    \coordinate (9d) at (9,4) {};
    \coordinate (9e) at (9,5) {};
    \coordinate (9f) at (9,6) {};
    \coordinate (9g) at (9,7) {};
    \coordinate (9h) at (9,8) {};
    \coordinate (9i) at (9,9) {};
    \coordinate (9j) at (9,10) {};
    \coordinate (10a) at (10,1) {};
    \coordinate (10b) at (10,2) {};
    \coordinate (10c) at (10,3) {};
    \coordinate (10d) at (10,4) {};
    \coordinate (10e) at (10,5) {};
    \coordinate (10f) at (10,6) {};
    \coordinate (10g) at (10,7) {};
    \coordinate (10h) at (10,8) {};
    \coordinate (10i) at (10,9) {};
    \coordinate (10j) at (10,10) {};
    
    \draw [black,thick] (1a) -- (1b);
    \draw [black,thick] (1b) -- (2b);
    \draw [black,thick] (2b) -- (2d);
    \draw [black,thick] (2d) -- (3d); 
    \draw [black,thick] (3d) -- (3e);
    \draw [black,thick] (3e) -- (4e);
    \draw [black,thick] (4e) -- (4f);
    \draw [black,thick] (4f) -- (6f);
    
    \draw [red,thick] (1a) -- (4a);
    \draw [red,thick] (4a) -- (4c);
    \draw [red,thick] (4c) -- (6c);
    \draw [red,thick] (6c) -- (6f);
    
    \draw [red,dotted,thick] (3a) -- (3c);
    \draw [red,dotted,thick] (3c) -- (4c);
       
   \end{tikzpicture}}
        \qquad
        \subfloat[]{\begin{tikzpicture}[scale = 0.5]
        \coordinate (1a) at (1,1) {};
    \coordinate (1b) at (1,2) {};
    \coordinate (1c) at (1,3) {};
    \coordinate (1d) at (1,4) {};
    \coordinate (1e) at (1,5) {};
    \coordinate (1f) at (1,6) {};
    \coordinate (1g) at (1,7) {};
    \coordinate (1h) at (1,8) {};
    \coordinate (1i) at (1,9) {};
    \coordinate (1j) at (1,10) {};
    \coordinate (2a) at (2,1) {};
    \coordinate (2b) at (2,2) {};
    \coordinate (2c) at (2,3) {};
    \coordinate (2d) at (2,4) {};
    \coordinate (2e) at (2,5) {};
    \coordinate (2f) at (2,6) {};
    \coordinate (2g) at (2,7) {};
    \coordinate (2h) at (2,8) {};
    \coordinate (2i) at (2,9) {};
    \coordinate (2j) at (2,10) {};
    \coordinate (3a) at (3,1) {};
    \coordinate (3b) at (3,2) {};
    \coordinate (3c) at (3,3) {};
    \coordinate (3d) at (3,4) {};
    \coordinate (3e) at (3,5) {};
    \coordinate (3f) at (3,6) {};
    \coordinate (3g) at (3,7) {};
    \coordinate (3h) at (3,8) {};
    \coordinate (3i) at (3,9) {};
    \coordinate (3j) at (3,10) {};
    \coordinate (4a) at (4,1) {};
    \coordinate (4b) at (4,2) {};
    \coordinate (4c) at (4,3) {};
    \coordinate (4d) at (4,4) {};
    \coordinate (4e) at (4,5) {};
    \coordinate (4f) at (4,6) {};
    \coordinate (4g) at (4,7) {};
    \coordinate (4h) at (4,8) {};
    \coordinate (4i) at (4,9) {};
    \coordinate (4j) at (4,10) {};
    \coordinate (5a) at (5,1) {};
    \coordinate (5b) at (5,2) {};
    \coordinate (5c) at (5,3) {};
    \coordinate (5d) at (5,4) {};
    \coordinate (5e) at (5,5) {};
    \coordinate (5f) at (5,6) {};
    \coordinate (5g) at (5,7) {};
    \coordinate (5h) at (5,8) {};
    \coordinate (5i) at (5,9) {};
    \coordinate (5j) at (5,10) {};
    \coordinate (6a) at (6,1) {};
    \coordinate (6b) at (6,2) {};
    \coordinate (6c) at (6,3) {};
    \coordinate (6d) at (6,4) {};
    \coordinate (6e) at (6,5) {};
    \coordinate (6f) at (6,6) {};
    \coordinate (6g) at (6,7) {};
    \coordinate (6h) at (6,8) {};
    \coordinate (6i) at (6,9) {};
    \coordinate (6j) at (6,10) {};
    \coordinate (7a) at (7,1) {};
    \coordinate (7b) at (7,2) {};
    \coordinate (7c) at (7,3) {};
    \coordinate (7d) at (7,4) {};
    \coordinate (7e) at (7,5) {};
    \coordinate (7f) at (7,6) {};
    \coordinate (7g) at (7,7) {};
    \coordinate (7h) at (7,8) {};
    \coordinate (7i) at (7,9) {};
    \coordinate (7j) at (7,10) {};
    \coordinate (8a) at (8,1) {};
    \coordinate (8b) at (8,2) {};
    \coordinate (8c) at (8,3) {};
    \coordinate (8d) at (8,4) {};
    \coordinate (8e) at (8,5) {};
    \coordinate (8f) at (8,6) {};
    \coordinate (8g) at (8,7) {};
    \coordinate (8h) at (8,8) {};
    \coordinate (8i) at (8,9) {};
    \coordinate (8j) at (8,10) {};
    \coordinate (9a) at (9,1) {};
    \coordinate (9b) at (9,2) {};
    \coordinate (9c) at (9,3) {};
    \coordinate (9d) at (9,4) {};
    \coordinate (9e) at (9,5) {};
    \coordinate (9f) at (9,6) {};
    \coordinate (9g) at (9,7) {};
    \coordinate (9h) at (9,8) {};
    \coordinate (9i) at (9,9) {};
    \coordinate (9j) at (9,10) {};
    \coordinate (10a) at (10,1) {};
    \coordinate (10b) at (10,2) {};
    \coordinate (10c) at (10,3) {};
    \coordinate (10d) at (10,4) {};
    \coordinate (10e) at (10,5) {};
    \coordinate (10f) at (10,6) {};
    \coordinate (10g) at (10,7) {};
    \coordinate (10h) at (10,8) {};
    \coordinate (10i) at (10,9) {};
    \coordinate (10j) at (10,10) {};
    
    \draw [black,thick] (1a) -- (1b);
    \draw [black,thick] (1b) -- (2b);
    \draw [black,thick] (2b) -- (2d);
    \draw [black,thick] (2d) -- (3d); 
    \draw [black,thick] (3d) -- (3e);
    \draw [black,thick] (3e) -- (4e);
    \draw [black,thick] (4e) -- (4f);
    \draw [black,thick] (4f) -- (6f);
    
    \draw [red,thick] (1a) -- (3a);
    \draw [red,thick] (3a) -- (3c);
    \draw [red,thick] (3c) -- (6c);
    \draw [red,thick] (6c) -- (6f);
    
    \draw [red,dotted,thick] (4c) -- (4d);
    \draw [red,dotted,thick] (4d) -- (6d);
       
   \end{tikzpicture}}
        \qquad
        \subfloat[]{\begin{tikzpicture}[scale = 0.5]
        \coordinate (1a) at (1,1) {};
    \coordinate (1b) at (1,2) {};
    \coordinate (1c) at (1,3) {};
    \coordinate (1d) at (1,4) {};
    \coordinate (1e) at (1,5) {};
    \coordinate (1f) at (1,6) {};
    \coordinate (1g) at (1,7) {};
    \coordinate (1h) at (1,8) {};
    \coordinate (1i) at (1,9) {};
    \coordinate (1j) at (1,10) {};
    \coordinate (2a) at (2,1) {};
    \coordinate (2b) at (2,2) {};
    \coordinate (2c) at (2,3) {};
    \coordinate (2d) at (2,4) {};
    \coordinate (2e) at (2,5) {};
    \coordinate (2f) at (2,6) {};
    \coordinate (2g) at (2,7) {};
    \coordinate (2h) at (2,8) {};
    \coordinate (2i) at (2,9) {};
    \coordinate (2j) at (2,10) {};
    \coordinate (3a) at (3,1) {};
    \coordinate (3b) at (3,2) {};
    \coordinate (3c) at (3,3) {};
    \coordinate (3d) at (3,4) {};
    \coordinate (3e) at (3,5) {};
    \coordinate (3f) at (3,6) {};
    \coordinate (3g) at (3,7) {};
    \coordinate (3h) at (3,8) {};
    \coordinate (3i) at (3,9) {};
    \coordinate (3j) at (3,10) {};
    \coordinate (4a) at (4,1) {};
    \coordinate (4b) at (4,2) {};
    \coordinate (4c) at (4,3) {};
    \coordinate (4d) at (4,4) {};
    \coordinate (4e) at (4,5) {};
    \coordinate (4f) at (4,6) {};
    \coordinate (4g) at (4,7) {};
    \coordinate (4h) at (4,8) {};
    \coordinate (4i) at (4,9) {};
    \coordinate (4j) at (4,10) {};
    \coordinate (5a) at (5,1) {};
    \coordinate (5b) at (5,2) {};
    \coordinate (5c) at (5,3) {};
    \coordinate (5d) at (5,4) {};
    \coordinate (5e) at (5,5) {};
    \coordinate (5f) at (5,6) {};
    \coordinate (5g) at (5,7) {};
    \coordinate (5h) at (5,8) {};
    \coordinate (5i) at (5,9) {};
    \coordinate (5j) at (5,10) {};
    \coordinate (6a) at (6,1) {};
    \coordinate (6b) at (6,2) {};
    \coordinate (6c) at (6,3) {};
    \coordinate (6d) at (6,4) {};
    \coordinate (6e) at (6,5) {};
    \coordinate (6f) at (6,6) {};
    \coordinate (6g) at (6,7) {};
    \coordinate (6h) at (6,8) {};
    \coordinate (6i) at (6,9) {};
    \coordinate (6j) at (6,10) {};
    \coordinate (7a) at (7,1) {};
    \coordinate (7b) at (7,2) {};
    \coordinate (7c) at (7,3) {};
    \coordinate (7d) at (7,4) {};
    \coordinate (7e) at (7,5) {};
    \coordinate (7f) at (7,6) {};
    \coordinate (7g) at (7,7) {};
    \coordinate (7h) at (7,8) {};
    \coordinate (7i) at (7,9) {};
    \coordinate (7j) at (7,10) {};
    \coordinate (8a) at (8,1) {};
    \coordinate (8b) at (8,2) {};
    \coordinate (8c) at (8,3) {};
    \coordinate (8d) at (8,4) {};
    \coordinate (8e) at (8,5) {};
    \coordinate (8f) at (8,6) {};
    \coordinate (8g) at (8,7) {};
    \coordinate (8h) at (8,8) {};
    \coordinate (8i) at (8,9) {};
    \coordinate (8j) at (8,10) {};
    \coordinate (9a) at (9,1) {};
    \coordinate (9b) at (9,2) {};
    \coordinate (9c) at (9,3) {};
    \coordinate (9d) at (9,4) {};
    \coordinate (9e) at (9,5) {};
    \coordinate (9f) at (9,6) {};
    \coordinate (9g) at (9,7) {};
    \coordinate (9h) at (9,8) {};
    \coordinate (9i) at (9,9) {};
    \coordinate (9j) at (9,10) {};
    \coordinate (10a) at (10,1) {};
    \coordinate (10b) at (10,2) {};
    \coordinate (10c) at (10,3) {};
    \coordinate (10d) at (10,4) {};
    \coordinate (10e) at (10,5) {};
    \coordinate (10f) at (10,6) {};
    \coordinate (10g) at (10,7) {};
    \coordinate (10h) at (10,8) {};
    \coordinate (10i) at (10,9) {};
    \coordinate (10j) at (10,10) {};
    
    \draw [black,thick] (1a) -- (1b);
    \draw [black,thick] (1b) -- (2b);
    \draw [black,thick] (2b) -- (2d);
    \draw [black,thick] (2d) -- (3d); 
    \draw [black,thick] (3d) -- (3e);
    \draw [black,thick] (3e) -- (4e);
    \draw [black,thick] (4e) -- (4f);
    \draw [black,thick] (4f) -- (6f);
    
    \draw [red,thick] (1a) -- (3a);
    \draw [red,thick] (3a) -- (3c);
    \draw [red,thick] (3c) -- (4c);
    \draw [red,thick] (4c) -- (4d);
    \draw [red,thick] (4d) -- (6d);
    \draw [red,thick] (6d) -- (6f);
    
    \draw [red,dotted,thick] (5d) -- (5e);
    \draw [red,dotted,thick] (5e) -- (6e);
       
   \end{tikzpicture}}
        \qquad
        \subfloat[]{\begin{tikzpicture}[scale = 0.5]
        \coordinate (1a) at (1,1) {};
    \coordinate (1b) at (1,2) {};
    \coordinate (1c) at (1,3) {};
    \coordinate (1d) at (1,4) {};
    \coordinate (1e) at (1,5) {};
    \coordinate (1f) at (1,6) {};
    \coordinate (1g) at (1,7) {};
    \coordinate (1h) at (1,8) {};
    \coordinate (1i) at (1,9) {};
    \coordinate (1j) at (1,10) {};
    \coordinate (2a) at (2,1) {};
    \coordinate (2b) at (2,2) {};
    \coordinate (2c) at (2,3) {};
    \coordinate (2d) at (2,4) {};
    \coordinate (2e) at (2,5) {};
    \coordinate (2f) at (2,6) {};
    \coordinate (2g) at (2,7) {};
    \coordinate (2h) at (2,8) {};
    \coordinate (2i) at (2,9) {};
    \coordinate (2j) at (2,10) {};
    \coordinate (3a) at (3,1) {};
    \coordinate (3b) at (3,2) {};
    \coordinate (3c) at (3,3) {};
    \coordinate (3d) at (3,4) {};
    \coordinate (3e) at (3,5) {};
    \coordinate (3f) at (3,6) {};
    \coordinate (3g) at (3,7) {};
    \coordinate (3h) at (3,8) {};
    \coordinate (3i) at (3,9) {};
    \coordinate (3j) at (3,10) {};
    \coordinate (4a) at (4,1) {};
    \coordinate (4b) at (4,2) {};
    \coordinate (4c) at (4,3) {};
    \coordinate (4d) at (4,4) {};
    \coordinate (4e) at (4,5) {};
    \coordinate (4f) at (4,6) {};
    \coordinate (4g) at (4,7) {};
    \coordinate (4h) at (4,8) {};
    \coordinate (4i) at (4,9) {};
    \coordinate (4j) at (4,10) {};
    \coordinate (5a) at (5,1) {};
    \coordinate (5b) at (5,2) {};
    \coordinate (5c) at (5,3) {};
    \coordinate (5d) at (5,4) {};
    \coordinate (5e) at (5,5) {};
    \coordinate (5f) at (5,6) {};
    \coordinate (5g) at (5,7) {};
    \coordinate (5h) at (5,8) {};
    \coordinate (5i) at (5,9) {};
    \coordinate (5j) at (5,10) {};
    \coordinate (6a) at (6,1) {};
    \coordinate (6b) at (6,2) {};
    \coordinate (6c) at (6,3) {};
    \coordinate (6d) at (6,4) {};
    \coordinate (6e) at (6,5) {};
    \coordinate (6f) at (6,6) {};
    \coordinate (6g) at (6,7) {};
    \coordinate (6h) at (6,8) {};
    \coordinate (6i) at (6,9) {};
    \coordinate (6j) at (6,10) {};
    \coordinate (7a) at (7,1) {};
    \coordinate (7b) at (7,2) {};
    \coordinate (7c) at (7,3) {};
    \coordinate (7d) at (7,4) {};
    \coordinate (7e) at (7,5) {};
    \coordinate (7f) at (7,6) {};
    \coordinate (7g) at (7,7) {};
    \coordinate (7h) at (7,8) {};
    \coordinate (7i) at (7,9) {};
    \coordinate (7j) at (7,10) {};
    \coordinate (8a) at (8,1) {};
    \coordinate (8b) at (8,2) {};
    \coordinate (8c) at (8,3) {};
    \coordinate (8d) at (8,4) {};
    \coordinate (8e) at (8,5) {};
    \coordinate (8f) at (8,6) {};
    \coordinate (8g) at (8,7) {};
    \coordinate (8h) at (8,8) {};
    \coordinate (8i) at (8,9) {};
    \coordinate (8j) at (8,10) {};
    \coordinate (9a) at (9,1) {};
    \coordinate (9b) at (9,2) {};
    \coordinate (9c) at (9,3) {};
    \coordinate (9d) at (9,4) {};
    \coordinate (9e) at (9,5) {};
    \coordinate (9f) at (9,6) {};
    \coordinate (9g) at (9,7) {};
    \coordinate (9h) at (9,8) {};
    \coordinate (9i) at (9,9) {};
    \coordinate (9j) at (9,10) {};
    \coordinate (10a) at (10,1) {};
    \coordinate (10b) at (10,2) {};
    \coordinate (10c) at (10,3) {};
    \coordinate (10d) at (10,4) {};
    \coordinate (10e) at (10,5) {};
    \coordinate (10f) at (10,6) {};
    \coordinate (10g) at (10,7) {};
    \coordinate (10h) at (10,8) {};
    \coordinate (10i) at (10,9) {};
    \coordinate (10j) at (10,10) {};
    
    \draw [black,thick] (1a) -- (1b);
    \draw [black,thick] (1b) -- (2b);
    \draw [black,thick] (2b) -- (2d);
    \draw [black,thick] (2d) -- (3d); 
    \draw [black,thick] (3d) -- (3e);
    \draw [black,thick] (3e) -- (4e);
    \draw [black,thick] (4e) -- (4f);
    \draw [black,thick] (4f) -- (6f);
    
    \draw [red,thick] (1a) -- (3a);
    \draw [red,thick] (3a) -- (3c);
    \draw [red,thick] (3c) -- (4c);
    \draw [red,thick] (4c) -- (4d);
    \draw [red,thick] (4d) -- (5d);
    \draw [red,thick] (6e) -- (6f);
    
    \draw [red,thick] (5d) -- (5e);
    \draw [red,thick] (5e) -- (6e);
       
   \end{tikzpicture}}
    \caption{Three consecutive applications of left-shifting on the red path. The figure on the far right is the furthest left we can shift the red path.}
    \label{fig:newleftshift1}
    \end{figure}

An example of left-shifting is given in Figure~\ref{fig:newleftshift1}. We are using left-shifting to modify the bottom path drawn in red. When we left-shift, we are not allowed to introduce a new intersection: hence how far we can shift a particular path to the left is bounded by other paths.


\subsection{Rectangular flip and reverse rectangular flip.}

\begin{Lemma} \label{lem:StartSplit} For any words $ B, C, D_1, \dots, D_m $ and an increasing sequence $A$,
\[
    \ell(AB,ABC,D_1, \dots, D_m) \leq \ell(AB,BC,D_1, \dots, D_m).
\]
\end{Lemma}

\begin{proof}
Let $S_{AB},S_{ABC},\ldots$ be a list of sequences that exhibits $\ell(AB,ABC,D_1, \dots, D_m)$. We will modify $ S_{AB},S_{ABC} $ into subsequences $  S_{AB}' ,  S_{ABC}' $ of $ ABC, B $ or $ AB, BC $, respectively, while maintaining their union, making the rest of subsequences irrelevant.

First, we partition $ S_{AB} $ and $ S_{ABC} $ into
\begin{align*}
    S_{AB} & = S_{AB}^1 S_{AB}^2, \\
    S_{ABC} & = S_{ABC}^1 S_{ABC}^2.
\end{align*}
where 
\[
    S_{AB}^1 \subseq A, \quad S_{AB}^2 \subseq B, \qquad S_{ABC}^1 \subseq A, \quad S_{ABC}^2 \subseq BC. 
\]
Let $S = S_{AB}^1 \sqcup S_{ABC}^1 $ listed in the same increasing order they appear in $ A $ so that $ S \subseq A $. Now, if $ \max(S_{AB}^1) < \max(S_{ABC}^1) $ or $ S_{AB}^1 $ is empty, let
\begin{align*}
    S_{ABC}' &= S \; S_{ABC}^2 \; \subseq ABC, \\
    S_{AB}' & = S_{AB}^2 \; \subseq B,
\end{align*}
and if $ \max(S_{ABC}^1) < \max(S_{AB}^1) $ or $ S_{ABC}^1 $ is empty, let
\begin{align*}
    S_{ABC}' &= S_{ABC}^2 \; \subseq BC, \\
    S_{AB}' & = S \; S_{AB}^2 \; \subseq AB,
\end{align*}
which are all increasing subsequences. Therefore, 
\begin{align*}
    \ell(AB,ABC,D_1, \dots, D_m) & \le \max( \ell(AB,BC,D_1, \dots, D_m), \ell(ABC,B,D_1, \dots, D_m)) \\
    & = \ell(AB,BC,D_1, \dots, D_m)
\end{align*}    
by \Cref{lem:LIS1}.
\end{proof}

We also have an analogous lemma where the proof is similar:

\begin{Lemma} \label{lem:EndSplit} For any words $ A, B, D_1, \dots, D_m $ and an increasing sequence $C$, 
\[
    \ell(ABC,BC,D_1, \dots, D_m) \leq \ell(AB,BC,D_1, \dots, D_m).
\]
\end{Lemma}

\begin{proof} The key idea in the proof of Lemma~\ref{lem:StartSplit} was that because $ A $ was increasing, we could move the elements of $A$ so that they are all contained in just one sequence instead of being split between $S_{AB}$ and $S_{ABC}$. The proof of \Cref{lem:EndSplit} works the same way as \Cref{lem:StartSplit} except we split according to the subsequences's restriction to $ C $ instead of $ A $ and then move every element of $ C $ to $S_{ABC}$ if $ \min(S_{ABC} \mid_C) < \min(S_{BC} \mid_C) $ or $ S_{BC} \mid_{C} $ is empty and to $S_{BC}$ otherwise.
\end{proof}

\begin{Lemma}[Rectangular flip] \label{lem:rectflip}
Suppose that the pair of paths $P,Q$ share the horizontal segment from $ (a,i) $ to $ (b, i) $. Suppose $ P $ continues to $(b,i+1)$ in the matrix and lies weakly above $ Q $. Then we may modify $ P $ to $ P' $ so that instead it goes from $(a,i)$ to $(a,i+1)$ to $(b,i+1)$, and then
\[
    \ell(P_1, \dots, P_k, P, Q) \le  \ell(P_1, \dots, P_k, P', Q).
\]
for any lattice paths $ P_1, \dots, P_k $.
\end{Lemma}

\begin{figure}[htbp]
    	\begin{tikzpicture}
\coordinate (1) at (-2,0) {};
\coordinate (2) at (1,0) {};
\coordinate (3) at (1,2) {};
\coordinate (4) at (2,2) {};
\coordinate (5) at (2,3) {};
\coordinate (6) at (3,3) {};
\coordinate (7) at (-2,-0.05) {};
\coordinate (8) at (1.05,-0.05) {};
\coordinate (9) at (1.05,0.95) {};
\coordinate (10) at (2.05,0.95){};
\coordinate (11) at (2.05,2.95){};
\coordinate (12) at (3,2.95){};
\draw [blue,thick](1) -- (2);
\draw [blue,thick](3) -- (2);
\draw [blue,thick](3) -- (4);
\draw [blue,thick](5) -- (4);
\draw [blue,thick](5) -- (6);

\node at (-2.5,0) {$(a, i)$};
\node at (1.5,0) {$(b, i)$};
\node at (0.3,1.2) {$(b, i + 1)$};

\node at (1,-0.3) {$A$};
\node at (1.3,0.5) {$B$};
\node at (1.3,1.5) {$C$};

\draw [fill = black] (1) circle (0.050);
\draw [fill = black] (2) circle (0.050);
\draw [fill = black] (1,1) circle (0.050);

\draw [red,thick](8) -- (9) ;
\draw [red,thick](9) -- (10) ;
\draw [red,thick](7) -- (8) ;
\draw [red,thick](10) -- (11) ;
\draw [red,thick](12) -- (11) ;
\coordinate (1') at (5,0) {};
\coordinate (2') at (5,1) {};
\coordinate (25) at (8,1) {};
\coordinate (3') at (8,2) {};
\coordinate (4') at (9,2) {};
\coordinate (5') at (9,3) {};
\coordinate (6') at (10,3) {};
\coordinate (7') at (5,-0.05) {};
\coordinate (8') at (8.05,-0.05) {};
\coordinate (9') at (8.05,0.95) {};
\coordinate (10') at (9.05,0.95){};
\coordinate (11') at (9.05,2.95){};
\coordinate (12') at (10,2.95){};
\draw [blue,thick](1') -- (2');
\draw [blue,thick](25) -- (2');
\draw [blue,thick](25) -- (3');
\draw [blue,thick](3') -- (4');
\draw [blue,thick](5') -- (4');
\draw [blue,thick](5') -- (6');

\node at (4.5,0) {$(a, i)$};
\node at (8.5,0) {$(b, i)$};
\node at (4.5,1.2) {$(a, i + 1)$};
\node at (7.3,1.2) {$(b, i + 1)$};

\node at (8,-0.3) {$A$};
\node at (8.3,0.5) {$B$};
\node at (8.3,1.5) {$C$};

\draw [fill = black] (1') circle (0.050);
\draw [fill = black] (2') circle (0.050);
\draw [fill = black] (8,0) circle (0.050);
\draw [fill = black] (8,1) circle (0.050);

\draw [red,thick](8') -- (9') ;
\draw [red,thick](9') -- (10') ;
\draw [red,thick](7') -- (8') ;
\draw [red,thick](10') -- (11') ;
\draw [red,thick](12') -- (11') ;

\end{tikzpicture}
\caption{A rectangular flip. Blue path: $ P \to P'$, Red path: $ Q $}
\label{fig:rectflip}
\end{figure}
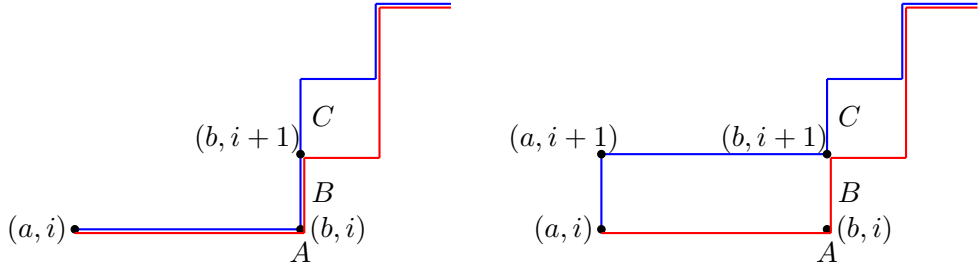

\begin{proof}
Let $ (R, S, \ldots) $ be disjoint increasing subsequences of $ P, Q, \ldots $ which exhibit the maximum combined length. Restricting to column $ b $ in both paths, let $ A = M_{b,i} $, $ B $ denote the portion of the column $ b $ where $ P $ and $ Q $ overlap, and $ C $ denote the portion of column $ b $ that only $ P $ contains, as in \Cref{fig:rectflip}. Note that $ B $ or $ C $ may be empty. Then, because $ A $ is increasing, we have
\[
    \ell(ABC, AB, \dots) \le \ell(BC, AB, \dots)
\]
by \Cref{lem:StartSplit} with $ ABC \subseq P $ and $ AB \subseq Q $. This means that we can assume without loss of generality that $ R \subseq P $ contains no elements from $ A = M_{b,i} $. Furthermore, since $ M_{a + 1,i} \ldots M_{b - 1,i} $ is increasing, we may move the elements of $ R \mid _{M_{a+1,i},\ldots,M_{b - 1,i}} $ to $ S $, without breaking that $ S $ is increasing. Together, we may assume without loss of generality that $ R $ uses none of the elements in $ M_{a + 1,i} \ldots M_{b,i} $. Hence, if we modify $ P $ to get $ P' $ by moving $ (a,i) - (b,i) - (b, i + 1) $ to $ (a,i) - (a,i + 1) - (b, i + 1) $ as in \Cref{fig:rectflip}, then $ P' $ still contains $ R $ as above. Thus, 
\[
    \ell(P_1, \dots, P_k, P, Q) \le  \ell(P_1, \dots, P_k, P', Q).
\]
\end{proof}

\begin{figure}
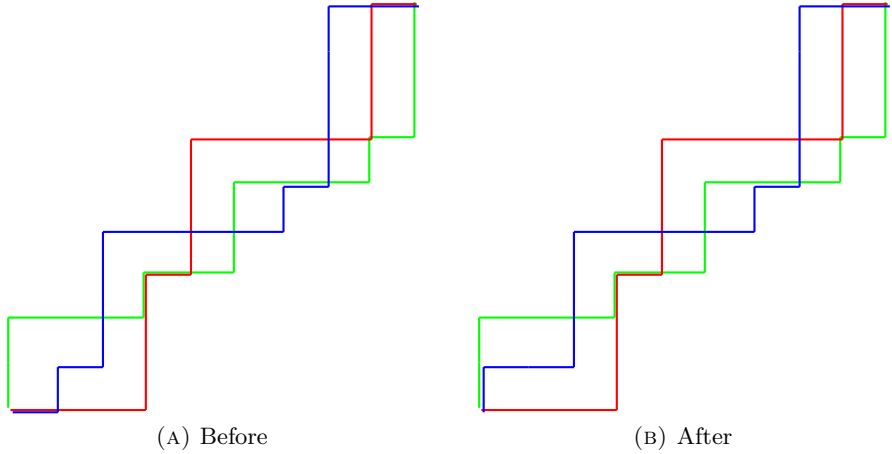

        \centering
        \subfloat[Before]{\input{Figures/newfig}}
        \qquad
        \subfloat[After]{\input{Figures/newrectflip}}
        
    \caption{An example of applying rectangular flip.}
    \label{fig:newrectflip}
    \end{figure}

An example of rectangular flip is given in Figure~\ref{fig:rectflip}. Note that the red and blue paths in the left figure have a common horizontal segment from $ (a,i) $ to $(b,i) $. After using the rectangular flip there, the blue path now is one height higher than its original position at the segment. An example of rectangular flip used in a family of paths is illustrated in Figure~\ref{fig:newrectflip}. Notice that in the left figure of Figure~\ref{fig:newrectflip}, the blue and red path share a horizontal segment at the beginning. After this operation, the green path goes through $(1,3)$, the blue goes through $(1,2)$ and the red path goes through $(1,1)$. In the next section, we will show at least one length-maximizing family of paths has the property of not just going through, but starting at, these points.

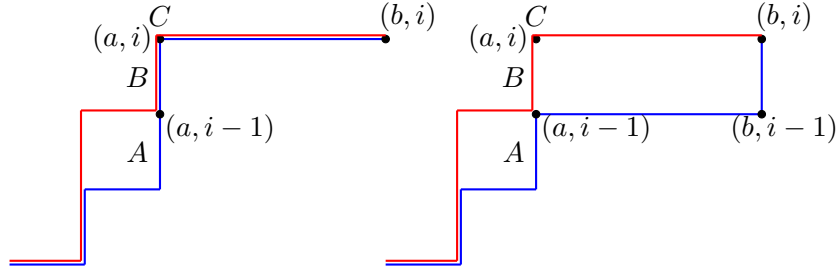
\begin{figure}[htbp]
    	\begin{tikzpicture}
\coordinate (1) at (0,10) {};
\coordinate (2) at (-3,10) {};
\coordinate (3) at (-3,8) {};
\coordinate (4) at (-4,8) {};
\coordinate (5) at (-4,7) {};
\coordinate (6) at (-5,7) {};
\coordinate (7) at (0,10.05) {};
\coordinate (8) at (-3.05,10.05) {};
\coordinate (9) at (-3.05,9.05) {};
\coordinate (10) at (-4.05,9.05){};
\coordinate (11) at (-4.05,7.05){};
\coordinate (12) at (-5,7.05){};
\draw [blue,thick](1) -- (2);
\draw [blue,thick](3) -- (2);
\draw [blue,thick](3) -- (4);
\draw [blue,thick](5) -- (4);
\draw [blue,thick](5) -- (6);

\node at (-3.5,10) {$(a,i)$};
\node at (-2.2,8.8) {$(a,i -1)$};
\node at (0.3,10.3) {$(b,i)$};

\node at (-3,10.3) {$C$};
\node at (-3.3,9.5) {$B$};
\node at (-3.3,8.5) {$A$};

\draw [fill = black] (1) circle (0.050);
\draw [fill = black] (2) circle (0.050);
\draw [fill = black] (-3,9) circle (0.050);

\draw [red,thick](8) -- (9) ;
\draw [red,thick](9) -- (10) ;
\draw [red,thick](7) -- (8) ;
\draw [red,thick](10) -- (11) ;
\draw [red,thick](12) -- (11) ;
\coordinate (1') at (5,10) {};
\coordinate (2') at (5,9) {};
\coordinate (25) at (2,9) {};
\coordinate (3') at (2,8) {};
\coordinate (4') at (1,8) {};
\coordinate (5') at (1,7) {};
\coordinate (6') at (0,7) {};
\coordinate (7') at (5,10.05) {};
\coordinate (8') at (1.95,10.05) {};
\coordinate (9') at (1.95,9.05) {};
\coordinate (10') at (0.95,9.05){};
\coordinate (11') at (0.95,7.05){};
\coordinate (12') at (0,7.05){};
\draw [blue,thick](1') -- (2');
\draw [blue,thick](25) -- (2');
\draw [blue,thick](25) -- (3');
\draw [blue,thick](3') -- (4');
\draw [blue,thick](5') -- (4');
\draw [blue,thick](5') -- (6');

\node at (1.5,10) {$(a,i)$};
\node at (2.8,8.8) {$(a, i - 1)$};
\node at (5.3,10.3) {$(b,i)$};
\node at (5.3,8.8) {$(b,i - 1)$};

\node at (2,10.3) {$C$};
\node at (1.7,9.5) {$B$};
\node at (1.7,8.5) {$A$};

\draw [fill = black] (1') circle (0.050);
\draw [fill = black] (2') circle (0.050);
\draw [fill = black] (2,9) circle (0.050);
\draw [fill = black] (2,10) circle (0.050);

\draw [red,thick](8') -- (9') ;
\draw [red,thick](9') -- (10') ;
\draw [red,thick](7') -- (8') ;
\draw [red,thick](10') -- (11') ;
\draw [red,thick](12') -- (11') ;

\end{tikzpicture}
\caption{A reverse rectangular flip. Blue path: $ P \to P'$, Red path: $ Q $}
\label{fig:revrectflip}
\end{figure}

\begin{Lemma}[Reverse rectangular flip] \label{lem:revrectflip}
Suppose that the pair of paths $P,Q$ share the horizontal segment from $ (a,i) $ to $ (b,i) $. Suppose $ P $ comes in via $(a,i-1)$ in the matrix and lies weakly below $ Q $. Then we may modify $ P $ to get $ P' $ that instead it goes from $(a,i-1)$ to $(b,i-1)$ to $(b,i)$, and then
\[
    \ell(P_1, \dots, P_k, P, Q) \le \ell(P_1, \dots, P_k, P', Q).
\]
for any lattice paths $ P_1, \dots, P_k $.
\end{Lemma}
\begin{proof} Let $ (R, S, \ldots) $ be disjoint increasing subsequences of $ P, Q, \ldots $ which exhibit the maximum possible length. Restricting to column $a$ in both paths, let $ C = M_{i,a} $, $ B $ denote the portion of the column $a$ where $ P $ and $ Q $ overlap, and $ A $ denote the portion of column $a$ that only $P$ contains as in \Cref{fig:revrectflip}. Note that $ A $ or $ B $ may be empty. Then, because $C$ is increasing, we have
\[
     \ell(ABC, BC, \dots) \le \ell(AB, BC, \dots)
\]
by \Cref{lem:EndSplit} with $ ABC \subseq P $ and $ BC \subseq Q $. Then, just as in the proof of \Cref{lem:rectflip}, we can assume without loss of generality that $ P' $ still contains $ R $. Hence, 
\[
    \ell(P_1, \dots, P_k, P, Q) \le \ell(P_1, \dots, P_k, P', Q).
\]

\end{proof}


\begin{figure}
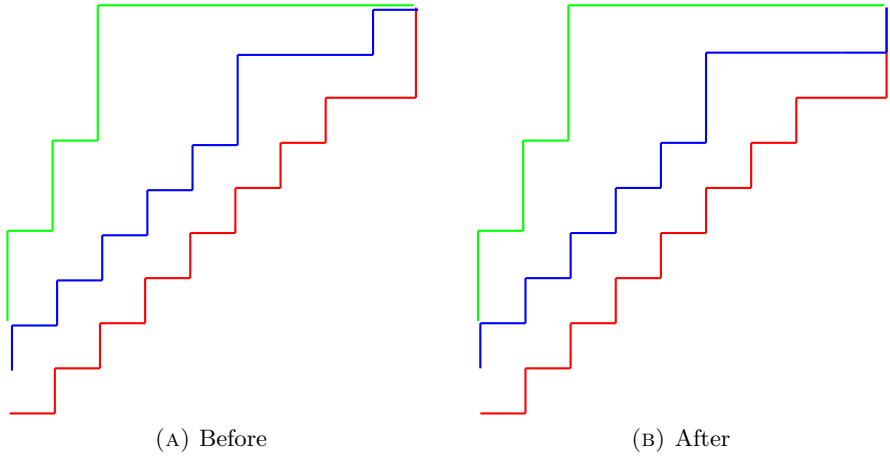

        \centering
        \subfloat[Before]{\input{Figures/newleftshift3}}
        \qquad
        \subfloat[After]{\input{Figures/newrectflip-1}}
        
    \caption{An example of applying reverse rectangular flip.}
    \label{fig:newrevrectflip}
    \end{figure}

An example of a reverse rectangular flip is given in Figure~\ref{fig:revrectflip}. Initially, the red and blue paths in the left figure have a common horizontal segment from $ (a,i) $ to $ (b,i) $. After using the reverse rectangular flip there, the blue path now is one height lower than its original position at the segment. An example of a reverse rectangular flip used in a family of paths is illustrated in Figure~\ref{fig:newrevrectflip}. Notice that in the left figure of Figure~\ref{fig:newrevrectflip}, the green and blue paths share a horizontal segment in the top row. After the reverse rectangular flip, the green and blue paths, the green path has a horizontal segment $(9,10)-(10,10)$, the blue has a horizontal segment $(9,9)-(10,9)$ and the red has $(9,8)-(10,8)$ which we will show in the next section is a property of at least one length-maximizing family of paths.

\section{The main result}
We have developed several tools regarding how the lattice paths encoding the longest sequences behave and how families of paths can be adjusted while weakly increasing the combined length they represent. Using these tools, we will be prove main result, \Cref{thm:main} in this section.

\begin{Lemma}
\label{lem:startingvertices}
There exists a family of $ k $ lattice paths $P_1,\ldots,P_k$ maximizing $ \ell(P_1, \dots, P_k) $ where for each $i$, the path $P_i$ starts at the point $(1,k+1-i)$ in the matrix.
\end{Lemma}
\begin{proof}
Let $P_1,\ldots,P_k$ be a family of lattice paths maximizing $ \ell(P_1, \dots, P_k) $. Add an auxiliary column 0 whose entries are empty words. Initially assume all paths start at (0,1) and go horizontally toward (1,1) without loss of generality. If not, extend them so they do. Order the paths based on their highest point in column 1, and label them $P_1,\ldots,P_k$ from top to bottom. We will perform rectangular flips to separate where these paths leave column 1. 

For $ j = 0, 1, \dots, (k - 1) $, say the paths $ P_1, \dots, P_k $ are in position $ j $ if paths $ P_1, \dots, P_{k - j} $ go from $ (0,1) $ to $ (0,j + 1) $ to $ (1, j + 1) $, and for each $ i \ge k - j $, $ P_i $ goes from $ (0,1) $ to $ (0, k - i  + 1) $ to $ (1, k - i + 1) $. Note both cases say the same thing for path $ P_{k - j} $. Initially, the paths start in position 0.

If the paths are in position $ j $, we will fix $ P_{k - j} $ and rectangular flip each of $ P_1, \dots, P_{k - j - 1} $ up out of the horizontal segment from $(0, j + 1)$ to $ (1, j + 1) $, which is already used by $ P_{k - j}$. After these rectangular flips, $ P_1, \dots, P_{k - j - 1} $ travel from $(0,0)$ to $ (0,j + 2) $ to $ (1,j + 2) $. Now our paths are in position $ j + 1 $.

Thus, we can inductively use rectangular flips to put the paths in position $ k - 1 $, where $ P_i $ travels from $ (0,1) $ to $ (0, k - i + 1) $ to $ (1, k - i + 1) $ for all $ i \ge 1 $. But column 0 is an auxiliary column that doesn't contain any letters to use, so we can remove it. When we remove column 0, the path $P_i$ starts at $(1,k - i + 1)$ for all $ i \ge 1 $, as desired.

\end{proof}

\begin{figure}[htbp]
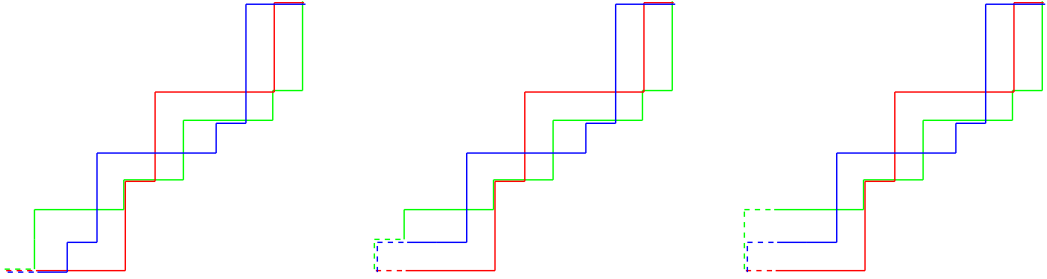

    \centering
     \scalebox{.66}{\input{Figures/newfigwithaux}}
     \qquad
     \scalebox{.66}{\input{Figures/newrectflipwithaux}}
     \qquad
     \scalebox{.66}{\input{Figures/newrectflipcopywithaux}}
     
     \caption{Running example: Start with a family of paths. Apply rectangular flip(s) to get the family satisfying Lemma~\ref{lem:startingvertices}. After this we may ignore the vertical segments in the first column to get a family satisfying the conditions of Lemma~\ref{lem:startingvertices}.}
     \label{fig:running1}
\end{figure}

Hence we may assume that each path $P_i$ starts at the point $(1,k - i + 1)$, as can be seen from the rightmost figure in Figure~\ref{fig:running1} (ignore the dashed lines). Beware that this does not necessarily mean each $P_i$ has the horizontal segment $(1,k - i + 1)-(2,k - i + 1)$.



\bs

We will reference a new family of paths - $ L_1, \dots, L_k $ that will serve as lower boundaries for our maximizing family. Fix positive integers $ q \ge r $, and $ k \ge 1 $. For each $ i = 1, \dots, k $, the \newword{bounding path} $ L_i $ is given by:

\begin{enumerate}[(a)]
\item Walk horizontally from $ (1, k-i + 1) $ to $ (i, k-i+1) $, then
\item Alternate between $ (0,1) $ and $ (1,0) $ steps, starting with $ (0,1) $, until reaching $ (r-k+i-1, r-i+1) $, then
 \item Walk horizontally to $ (q,r-i+1) $.
\end{enumerate}

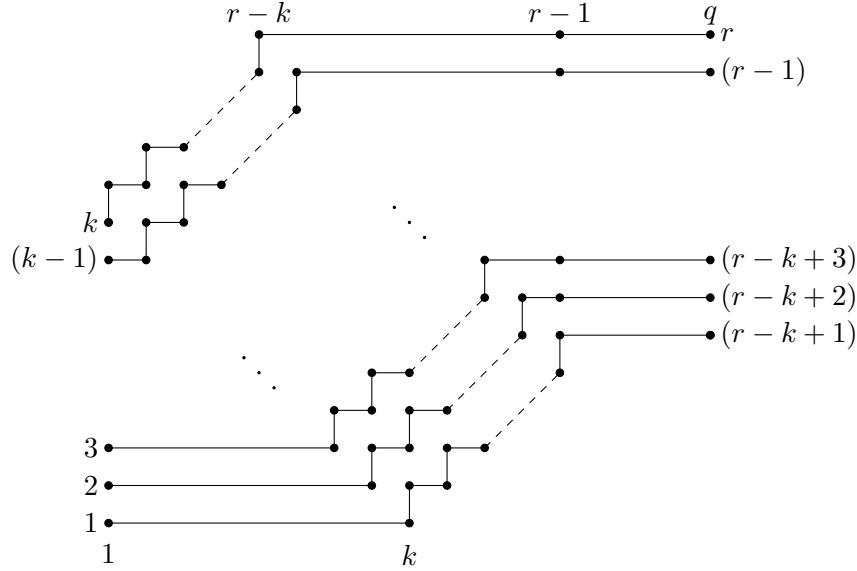
\begin{figure} 
\begin{tikzpicture}[scale=1, transform shape]

\coordinate [label=left:$1$] (1a) at (0,0) {};
\coordinate (1b) at (4,0) {};
\coordinate [label=below:$k$] (1b') at (4,-0.15) {};
\coordinate (1c) at (4,0.5) {};
\coordinate (1d) at (4.5,0.5) {};
\coordinate (1e) at (4.5,1) {};
\coordinate (1f) at (5,1) {};
\coordinate (1g) at (6,2) {};
\coordinate (1h) at (6,2.5) {};
\coordinate (1i) at (8,2.5) {};
\coordinate (1j) at (8,2.5) {};

\coordinate [label=left:$2$] (2a) at (0,0.5) {};
\coordinate (2b) at (3.5,0.5) {};
\coordinate (2c) at (3.5,1) {};
\coordinate (2d) at (4,1) {};
\coordinate (2e) at (4,1.5) {};
\coordinate (2f) at (4.5,1.5) {};
\coordinate (2g) at (5.5, 2.5) {};
\coordinate (2h) at (5.5, 3) {};
\coordinate (2i) at (6, 3) {};
\coordinate (2j) at (8, 3) {};

\coordinate [label=left:$3$] (3a) at (0,1) {};
\coordinate (3b) at (3,1) {};
\coordinate (3c) at (3,1.5) {};
\coordinate (3d) at (3.5,1.5) {};
\coordinate (3e) at (3.5,2) {};
\coordinate (3f) at (4,2) {};
\coordinate (3g) at (5, 3) {};
\coordinate (3h) at (5, 3.5) {};
\coordinate (3i) at (6, 3.5) {};
\coordinate (3j) at (8,3.5) {};

\coordinate [label=left:$(k - 1)$] (4a) at (0,3.5) {};
\coordinate (4b) at (0.5,3.5) {};
\coordinate (4c) at (0.5, 4) {};
\coordinate (4d) at (1, 4) {};
\coordinate (4e) at (1,4.5) {};
\coordinate (4f) at (1.5, 4.5) {};
\coordinate (4g) at (2.5, 5.5) {};
\coordinate (4h) at (2.5, 6) {};
\coordinate (4i) at (6, 6) {};
\coordinate (4j) at (8, 6) {};

\coordinate [label=left:$k$] (5a) at (0,4) {};
\coordinate (5b) at (0,4) {};
\coordinate (5c) at (0,4.5) {};
\coordinate (5d) at (0.5, 4.5) {};
\coordinate (5e) at (0.5,5) {};
\coordinate (5f) at (1, 5) {};
\coordinate (5g) at (2, 6) {};
\coordinate [label=above:$r - k$] (5h) at (2, 6.5) {};
\coordinate [label=above:$r - 1$] (5i) at (6, 6.5) {};
\coordinate [label=above:$q$] (5j) at (8, 6.5) {};

\coordinate [label=right: $ r $] (6) at (8,6.5) {};
\coordinate [label=below: $1$] (7) at (0,-0.15) {};
\coordinate [label=right: $ (r - k + 1)$] (8) at (8,2.5) {};
\coordinate [label=right: $ (r - k + 2)$] (9) at (8,3) {};
\coordinate [label=right: $ (r - k + 3)$] (11) at (8,3.5) {};
\coordinate [label=right: $ (r - 1) $] (10) at (8,6) {};

\draw (1a) -- (1b);
\draw (1b) -- (1c);
\draw (1c) -- (1d);
\draw (1d) -- (1e);
\draw (1e) -- (1f);
\draw[style = dashed] (1f) -- (1g);
\draw (1g) -- (1h);
\draw (1h) -- (1i);
\draw (1i) -- (1j);

\draw (2a) -- (2b);
\draw (2b) -- (2c);
\draw (2c) -- (2d);
\draw (2d) -- (2e);
\draw (2e) -- (2f);
\draw[style = dashed] (2f) -- (2g);
\draw (2g) -- (2h);
\draw (2h) -- (2i);
\draw (2i) -- (2j);

\draw (3a) -- (3b);
\draw (3b) -- (3c);
\draw (3c) -- (3d);
\draw (3d) -- (3e);
\draw (3e) -- (3f);
\draw[style = dashed] (3f) -- (3g);
\draw (3g) -- (3h);
\draw (3h) -- (3i);
\draw (3i) -- (3j);

\draw (4a) -- (4b);
\draw (4b) -- (4c);
\draw (4c) -- (4d);
\draw (4d) -- (4e);
\draw (4e) -- (4f);
\draw[style = dashed] (4f) -- (4g);
\draw (4g) -- (4h);
\draw (4h) -- (4i);
\draw (4i) -- (4j);

\draw (5a) -- (5b);
\draw (5b) -- (5c);
\draw (5c) -- (5d);
\draw (5d) -- (5e);
\draw (5e) -- (5f);
\draw[style = dashed] (5f) -- (5g);
\draw (5g) -- (5h);
\draw (5h) -- (5i);
\draw (5i) -- (5j);

\draw [fill = black] (1a) circle (0.050);
\draw [fill = black] (1b) circle (0.050);
\draw [fill = black] (1c) circle (0.050);
\draw [fill = black] (1d) circle (0.050);
\draw [fill = black] (1e) circle (0.050);
\draw [fill = black] (1f) circle (0.050);
\draw [fill = black] (1g) circle (0.050);
\draw [fill = black] (1h) circle (0.050);
\draw [fill = black] (1i) circle (0.050);
\draw [fill = black] (1j) circle (0.050);

\draw [fill = black] (2a) circle (0.050);
\draw [fill = black] (2b) circle (0.050);
\draw [fill = black] (2c) circle (0.050);
\draw [fill = black] (2d) circle (0.050);
\draw [fill = black] (2e) circle (0.050);
\draw [fill = black] (2f) circle (0.050);
\draw [fill = black] (2g) circle (0.050);
\draw [fill = black] (2h) circle (0.050);
\draw [fill = black] (2i) circle (0.050);
\draw [fill = black] (2j) circle (0.050);

\draw [fill = black] (3a) circle (0.050);
\draw [fill = black] (3b) circle (0.050);
\draw [fill = black] (3c) circle (0.050);
\draw [fill = black] (3d) circle (0.050);
\draw [fill = black] (3e) circle (0.050);
\draw [fill = black] (3f) circle (0.050);
\draw [fill = black] (3g) circle (0.050);
\draw [fill = black] (3h) circle (0.050);
\draw [fill = black] (3i) circle (0.050);
\draw [fill = black] (3j) circle (0.050);

\draw [fill = black] (4a) circle (0.050);
\draw [fill = black] (4b) circle (0.050);
\draw [fill = black] (4c) circle (0.050);
\draw [fill = black] (4d) circle (0.050);
\draw [fill = black] (4e) circle (0.050);
\draw [fill = black] (4f) circle (0.050);
\draw [fill = black] (4g) circle (0.050);
\draw [fill = black] (4h) circle (0.050);
\draw [fill = black] (4i) circle (0.050);
\draw [fill = black] (4j) circle (0.050);

\draw [fill = black] (5a) circle (0.050);
\draw [fill = black] (5b) circle (0.050);
\draw [fill = black] (5c) circle (0.050);
\draw [fill = black] (5d) circle (0.050);
\draw [fill = black] (5e) circle (0.050);
\draw [fill = black] (5f) circle (0.050);
\draw [fill = black] (5g) circle (0.050);
\draw [fill = black] (5h) circle (0.050);
\draw [fill = black] (5i) circle (0.050);
\draw [fill = black] (5j) circle (0.050);

\draw [fill = black] (2,2) circle (0.015);
\draw [fill = black] (2.2,1.8) circle (0.015);
\draw [fill = black] (1.8,2.2) circle (0.015);
\draw [fill = black] (4,4) circle (0.015);
\draw [fill = black] (4.2,3.8) circle (0.015);
\draw [fill = black] (3.8,4.2) circle (0.015);

\end{tikzpicture}
\caption{The lower boundary paths for $ k $ paths: $ L_1, \dots, L_k$ (from top to bottom). The numbers on the left and right indicate the row that the corresponding points are on. The numbers on the top and bottom indicate the columns the corresponding points are on.} 
\label{fig:lowerboundpaths}
\end{figure}

\Cref{fig:lowerboundpaths} shows $ L_1, \dots, L_k $. These paths appear in \cite[Lemma 4.2]{TabStab} as well. If all rows of $ T $ have the same size, this family of paths is a length-maximizing family in that 
\[
    \ell(P_1, \dots, P_k) \le \ell(L_1,\dots, L_k)
\]
for any paths $ P_1, \dots P_k $ \cite[Lemma 4.2]{TabStab}. If the rows of $ T $ have different sizes, this may no longer be the case, and we have to allow our family of paths to be weakly above these in respective top to bottom order.

\begin{Lemma}
\label{lem:horizseg}
There exists a family of paths $P_1,\ldots,P_k$ that maximizes $ \ell(P_1, \dots, P_k) $ such that $ P_i $ is weakly above $ L_i $ for all $ i $, and each $P_i$ contains the horizontal segment $(r - k + i - 1,r-i+1) - (q,r-i+1)$.
\end{Lemma}
\begin{proof}
 Let $P_1,\ldots,P_k$ be a family of lattice paths maximizing $ \ell(P_1, \dots, P_k) $. From Lemma~\ref{lem:startingvertices} we may assume that each $P_i$ starts at $(1,k-i+1)$. From repeated usage of top-down switching lemma we may assume that all paths do not cross (they may intersect) and hence come in a top to bottom order, say $ P_1, \dots, P_k $ from top to bottom.  

We claim that for each $i$, we can make $P_i$ stays weakly above the bounding path $L_i$. For path $P_1$, left-shifting gives us some path that starts at $(1,k)$ with an up move, and each horizontal segment has size $1$ before the path reaches the top row. Hence $P_1$ is weakly above $ L_1 $. Now inductively assume that path $ P_j $ lies above $ L_j $ for all $ j = 1, \dots, (i - 1) $. If $P_i$ dips below $L_i$ at some point, pick the earliest point that this happens. Then in this row, $P_i$ must have horizontal segment of size at least $2$. From the fact that $P_{i-1}$ is weakly above $L_{i-1}$, we have the empty space necessary to left-shift $P_i$ weakly above $L_i$ in this row using \Cref{lem:leftshift}. We continue to similarly adjust $ P_i $ by left-shifting any other points where $ P_i $ goes below $ L_i $, so we may assume without loss of generality that $ P_i $ stays weakly above $ L_i $, completing our induction.

The last step is to show that each path $P_i$ has the horizontal segment $(r - k + i - 1,r-i+1) - (q,r -i + 1)$. We again employ induction on $i$. For path $ P_1 $, which stays weakly above $L_1$, it has the desired horizontal segment because $ L_1 $ has it in the top row. Inductively assume that $P_{i-1}$ has the horizontal segment $(r - k + i - 2, r - i + 2) - (q,r - i + 2)$. As $ P_i $ is weakly below $ P_{i - 1} $ but weakly above $ L_i $, $ P_i $ must pass through $ (r - k + i - 1,r -i + 1) $ and end with a horizontal segment in either row $ r - i + 2 $ or row $ r -i + 1 $. Now for any of $ P_i, P_{i + 1}, \dots, P_k $ that end in row $ r - i + 2 $, we can reverse rectangular flip it, \Cref{lem:revrectflip}, so it ends in row $ r - i + 1 $ instead because $ P_{i - 1} $ already uses the same portion of row $ r - i + 2 $. This maintains the top-down order of the paths, so we can and will use the same labels for the paths. Thus, we may assume without loss of generality that $ P_i $ ends in row $ r -i + 1 $. As it is also weakly above $ L_i $, $P_i$ must end in the horizontal segment $(r - k + i - 1,r -i + 1) - (q, r -i + 1)$, completing our induction. 

\end{proof}

\begin{figure}[htbp]
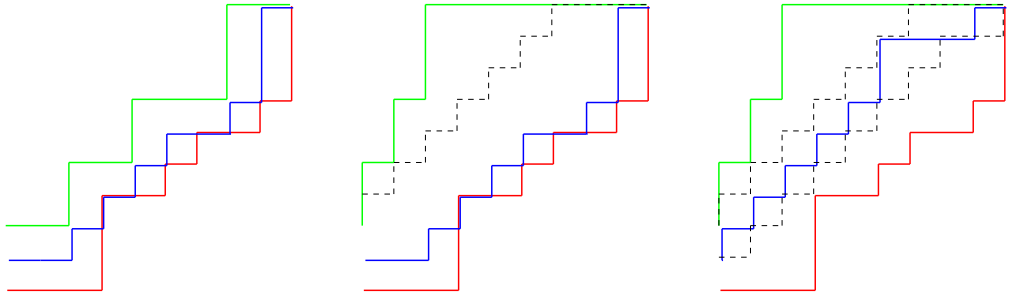

    \centering
     \scalebox{.7}{\input{Figures/newtopdown}}
     \qquad
     \scalebox{.7}{\input{Figures/newleftshift1}}
     \qquad
     \scalebox{.7}{\input{Figures/newleftshift2}}
     \caption{Running example: Apply top-down switching to get a family of paths that do not cross. Next apply left-shift on $P_1$ to get a path that stays weakly above the bounding path $L_1$, drawn in black dotted lines. Do the same for $P_2$ afterwards.}
     \label{fig:running2}
\end{figure}

\begin{figure}[htbp]
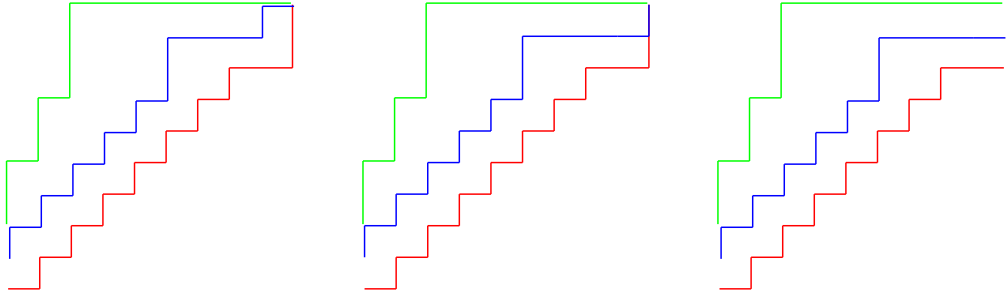

     \centering
     \scalebox{0.7}{\input{Figures/newleftshift3}}
     \qquad
     \scalebox{0.7}{\input{Figures/newrectflip-1}}
     \qquad
     \scalebox{0.7}{\input{Figures/newrectflip-2}}
\caption{Running example: Apply left-shift on $P_3$. Lastly, apply reverse rectangular flip(s) to get the family satisfying Lemma~\ref{lem:horizseg}.}
\label{fig:running3}
\end{figure}

For example, start with the family of paths in the left figure of Figure~\ref{fig:running1}. Using rectangular flips as in Lemma~\ref{lem:startingvertices} in an auxiliary column gives us the figure in the middle and then the figure on the right, where the paths satisfy the property of Lemma~\ref{lem:startingvertices}. The dashed lines are part of the auxiliary column, which play a part in the rectangular flips. As they do not contribute anything to the increasing subsequences, we remove them afterward. Hence we may assume the paths start at $(1,1),(1,2),(1,3)$ as in the rightmost figure of Figure~\ref{fig:running1}, ignoring the dashed lines.

Next using top-down switching lemma here gives us a family where we can label the paths $P_1,P_2,P_3$ from top to bottom in a way so that path $P_i$ is weakly above path $P_{i+1}$ for each $i$ as in the leftmost figure of Figure~\ref{fig:running2}. Next use the left-shifting on the top path $P_1$ to get the middle figure of Figure~\ref{fig:running2}. Notice the resulting modification of $P_1$ gives us a path that lies weakly above the bounding path $L_1$ drawn in black dashed lines. Next use left-shifting on $P_2$ to get the rightmost figure, where the resulting modification of $P_2$ lies weakly above the bounding path $L_2$. Now do the same for $P_3$ as well to get the left figure of Figure~\ref{fig:running3}. Lastly use reverse rectangular flips at the end to get the desired family of Lemma~\ref{lem:horizseg}.


Now we are ready to prove our main result.

\begin{customthm}{\ref{thm:main}}
For any skew standard tableaux $ T $ with $ r $ rows and weakly decreasing row sizes from top to bottom,
\[
    \stab(T) \le r.
\]
\end{customthm}
\begin{proof}

Let $ \lam_1 \ge \dots \ge \lam_r $ be the row sizes of $ T $ from top to bottom. By Greene's theorem, \cite{MR0354395}, we have that
\[
    \sum_{j = 1}^k \sh(\Rect(T^{(q)}))_j = \ell( \underbrace{w(q,T), \dots, w(q,T)}_{k \tx{ times} } ).
\]
Recall $ w(q,T) $ represents the reading word of $ T^{(q)} $.

Now, for each $ q \ge r - 1 $, \Cref{lem:horizseg} tells us that $ \ell( \underbrace{w(q,T), \dots, w(q,T)}_{k \tx{ times} } ) $ is exhibited by increasing subsequences of some family of lattice paths $P_1, \dots, P_k $ where $ P_i $ ends in the horizontal segment $(r - k + i - 1,r-i+1) - (q,r-i+1)$ for all $ i = 1, \dots, k $. Note that because $ i \le k $, we always have $ r - k + i - 1 \le r - 1 \le q $.

In particular, to go from $ q = r - 1 $ to $ q = r $, we may simply add on the increasing sequences represented by the points $ (r,r), (r, r - 1), \dots, (r, r - k + 1) $ to a length-maximizing family for $ q = r  - 1 $ from top to bottom. Note these points represent sequences with sizes $ \lam_1, \dots, \lam_k $ respectively. This family of disjoint increasing sequences of $ w(r,T) $ tells us
\begin{equation}
\begin{aligned}
\label{eq:r-1tor}
    \sum_{j = 1}^k \sh(\Rect(T^{(r)}))_j & = \ell(\underbrace{w(r,T), \dots, w(r,T)}_{k \tx{ times} } ) \\
    & \ge \ell(\underbrace{w(r - 1,T), \dots, w(r-1,T)}_{k \tx{ times} } ) + \lam_1 + \dots + \lam_k \\
    & = \sum_{j = 1}^k \sh(\Rect(T^{(r - 1)}))_j + \lam_1 + \dots + \lam_k.
\end{aligned}
\end{equation}

On the other hand, if we remove the the points $ (r,r), (r, r - 1), \dots, (r, r - k + 1) $ from such a length-maximizing family of lattice paths for $ q = r $ - we can assume the family contains these points by \Cref{lem:horizseg} - we get a family of paths for $ q = r - 1 $. Note that we removed $ \lam_1 + \dots + \lam_k $ elements potentially used in increasing sequences. In terms of increasing subsequences, this means 
\begin{equation}
\begin{aligned}
\label{eq:rtor-1}
    \sum_{j = 1}^k \sh(\Rect(T^{(r - 1)}))_j & = \ell(\underbrace{w(r - 1,T), \dots, w(r-1,T)}_{k \tx{ times} }) \\
    & \ge \ell(\underbrace{w(r,T), \dots, w(r,T)}_{k \tx{ times} } ) - (\lam_1 + \dots + \lam_k) \\
    & = \sum_{j = 1}^k \sh(\Rect(T^{(r)}))_j - (\lam_1 + \dots + \lam_k).
\end{aligned}
\end{equation}
Putting \eqref{eq:r-1tor} and \eqref{eq:rtor-1} together, 
\begin{align}
\label{eq:shapeincrease}
    \sum_{j = 1}^k \sh(\Rect(T^{(r)}))_j - \sum_{j = 1}^k \sh(\Rect(T^{(r - 1)}))_j = \lam_1 + \dots + \lam_k.
\end{align}
Since \eqref{eq:shapeincrease} holds for all $ k = 1, \dots, r $, we must have
\[
     \sh(\Rect(T^{(r)}))_j - \sh(\Rect(T^{(r - 1)}))_j = \lam_j
\]
for all $ j = 1, \dots, r $. Hence the rows of $ \Rect(T^{(r)}) \mid_{[(r - 1)n + 1, rn]} $ have the same sizes as the rows of $ T + (r - 1)n $, namely $ \lam_1, \lam_2, \dots, \lam_r $ from top to bottom. But this can only happen if every element of $ T + (r - 1)n $ experiences no vertical slides during rectification. Therefore, $ T $ stabilizes at $ r $, so $ \stab(T) \le r $.
\end{proof}

The main idea of our proof was to narrow our search for a length-maximizing family of lattice paths to those satisfying a certain property - namely lying above the bounding paths and ending in the same horizontal segments at the bounding paths. Is there is an efficient way to construct a particular family of length-maximizing paths given the tableau?

\begin{Problem}
Given a skew tableau $T$, is there a simple method to find a length-maximizing family of $k$-paths in $M(q,T)$ for each $k$ and $q$?
\end{Problem}

\begin{figure}
   \input{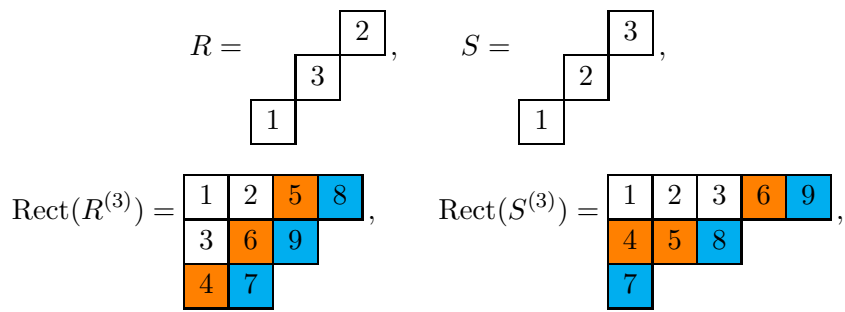}
   \caption{Tableaux $R=T_{132}$ and $S=T_{123}$. We can see that the stabilization index of $132$ is $2$ whereas the stabilization index of $123$ is $3$.}
\label{fig:permutationexample}
\end{figure}

We have found an upper bound on the stabilization index of a tableau. It seems a far-fetched goal at this moment to obtain this index without using Greene's theorem or constructing the rectification. A reasonable subclass of tableaux to see if such thing is possible would be tableaux constructed from permutations as in \Cref{fig:permutationexample}. Given a permutation $w$, we can construct the skew tableau $T_w$ which has one entry per each row and $w$ as its reading word. Then we can define $\stab(w)$ as $\stab(T_w)$, which give us a permutation statistic!

In Chapter 8 of \cite{TabStab}, Ahlbach introduced the stabilization index as a permutation statistic. He showed that $ \stab(w) $ is bounded strictly below by the ascent statistic, \cite[Lemma 8.4]{TabStab}, showed that $ \stab(w) $ depends only the recording tableau $ Q(w) $, \cite[Lemma 8.3]{TabStab}, characterized the permutations with stab 1, \cite[Lemma 8.5]{TabStab}, and characterized the permutations with stab 2, \cite[Theorem 8.7]{TabStab}. But there are still many open questions! A full characterization of $ \stab(w) $ in terms of $ Q(w) $ would be ideal.

\begin{Problem}[\cite{TabStab}]
 For a permutation $ w $, is there a way to find $ \stab(w) $ directly from the permutation or its recording tableau (that is, without constructing the rectification of or using Greene's theorem)? What is the relationship between $ \stab(w) $ and $ Q(w) $?
\end{Problem}

\section*{Acknowledgements}
This research was primarily conducted during the 2020 Honors Summer Math Camp at Texas State University. The authors gratefully acknowledge the support from the camp and also thank Texas State University for providing support and a great working environment. We would also like to thank the reviewer for their helpful comments and for pointing out relevant references.

\bibliographystyle{plain}    
\bibliography{refs}

\end{document}